\documentclass[12pt]{article}
\usepackage{amsmath}
\usepackage{amssymb}
\usepackage{url}
\newtheorem{thm}{Theorem}
\newtheorem{lemma}[thm]{Lemma}
\newtheorem{cor}[thm]{Corollary}
\newtheorem{prop}[thm]{Proposition}
\newtheorem{definition}[thm]{Definition}

\DeclareMathOperator{\esssup}{ess\,sup}

\def\vint_#1{\mathchoice%
          {\mathop{\kern 0.2em\vrule width 0.6em height 0.69678ex depth -0.58065ex
                  \kern -0.8em \intop}\nolimits_{\kern -0.4em#1}}%
          {\mathop{\kern 0.1em\vrule width 0.5em height 0.69678ex depth -0.60387ex
                  \kern -0.6em \intop}\nolimits_{#1}}%
          {\mathop{\kern 0.1em\vrule width 0.5em height 0.69678ex
              depth -0.60387ex
                  \kern -0.6em \intop}\nolimits_{#1}}%
          {\mathop{\kern 0.1em\vrule width 0.5em height 0.69678ex depth -0.60387ex
                  \kern -0.6em \intop}\nolimits_{#1}}}

\newcommand{\Rn}{\ensuremath{\mathbb{R}^{n}}}
\newcommand{\R}{\ensuremath{\mathbb{R}}}

\begin{document}
\title{\textsc{Notes on the\\ Infinity-Laplace Equation}}
\author{Peter Lindqvist}
\date{{\small Norwegian University of Science and Technology}}
\maketitle

\bigskip

{\footnotesize \textsf{These notes are written up after my lectures at the University of Pittsburgh in March 2014 and at
      Tsinghua University  in May 2014. My objective is the $\infty$-Laplace Equation, a marvellous kin to the ordinary Laplace Equation. The $\infty$-Laplace Equation  has delightful counterparts to the Dirichlet integral, the Mean Value Theorem, the Brownian Motion, Harnack's Inequality and so on. It has applications to image processing and to mass transfer problems and provides optimal Lipschitz extensions of boundary values. My treaty of this ''fully non-linear'' degenerate equation is far from complete and  generalizations are deliberately avoided.  ---''Less is more.''}}

\bigskip

\qquad\quad\, \texttt{Habent sua fata libelli pro captu lectoris.}

\section{Introduction}

These notes are about the celebrated $\infty$-Laplace Equation, which so far as I know first was derived by G. Aronsson in 1967, see [A1]. This fascinating equation has later been rediscovered in connexion with Image Processing and in Game Theory. In two variables it reads
$$u_x^2u_{xx}+ 2u_xu_yu_{xy}+u_{y}^2u_{yy}\,=\,0,$$
where $u = u(x,y).$\footnote{Formally it would be classified as parabolic, since its discriminant is $u_x^2u_y^2 -(u_xu_y)^2 = 0,$ but the main features are ''elliptic''. The characteristic curves are $u_xdy-u_ydx = 0.$ In several variables, it is degenerate elliptic.} For a function $u=u(x_1,x_2,\dots,x_n)$ in $n$  variables the equation is$$
\boxed{\Delta_{\infty}u \,\equiv \,\sum_{i,j=1}^{n}\frac{\partial u}{\partial x_i} \frac{\partial u}{\partial x_j}\frac{\partial^2u}{\partial x_i \partial x_j}\,=\,0.}$$ It is the limit equation of the $p$-Laplace Equations 
\begin{equation}
\label{pl}
\Delta_pu \equiv |\nabla u|^{p-4}\left\{|\nabla u|^2\Delta u + (p-2)\Delta_{\infty}u\right\}\,=\,0 
\end{equation}
as $p \to \infty.$ To see this, divide out the factor $|\nabla u|^{p-4},$ divide by $p-2$ and let $p \to \infty$ in the resulting expression
$$\dfrac{|\nabla u|^2\Delta u}{p-2} + \Delta_{\infty}u \,=\,0.$$
Our derivation leaves much to be desired, but the limit equation is the correct one. ---The solutions are often called $\infty$-harmonic functions.

\paragraph{Classical solutions}  A function $u \in C^2(\Omega),\,\,\Omega$ denoting a domain in $n$-dimensional Euclidean space $\Rn,$ which satisfies the equation pointwise is called a classical solution. A few explicit examples are 
\begin{align*}
b|x-\mathbf{a}| + c,\qquad    a\sqrt{x_1^2+x_2^2 +\cdots+x_k^2} + b,\qquad \langle \mathbf{a},x\rangle + b,\\
\arctan\frac{x_2}{x_1},\qquad \arctan \frac{x_3}{\sqrt{x_1^2+x^2_2}},\qquad x_1^{\frac{4}{3}}-x_2^{\frac{4}{3}}, \quad \text {and} \quad\\
a_1x_1^{\frac{4}{3}} + a_2x_2^{\frac{4}{3}} + \cdots + a_n x_n^{\frac{4}{3}}\quad \text{where}\quad a_1^3+a_2^3+\cdots+ a_n^3 = 0,
\end{align*}
in domains where they are of class $C^2.$ Since
$$\boxed{\Delta_{\infty}u = \frac{1}{2}\langle \nabla u, \nabla |\nabla u|^2 \rangle}$$
all solutions $u \in C^1$ of the \emph{eikonal equation} $|\nabla u|^2 = \mathrm{C}$ can be included. Finally, solutions in \emph{disjoint variables} can be superposed:
$$\sqrt{x^2_1+x^2_2}-7\sqrt{x^2_3+x^2_4} + x_5 +x_6^{\frac{4}{3}} -x_7^{\frac{4}{3}},$$
for example. Also solutions like
$$|x_{i_{1}}^2+\cdots x_{i_{\nu}}^2|^{\frac{2}{3}}-|x_{j_{1}}^2+\cdots x_{j_{\mu}}^2|^{\frac{2}{3}} $$
are possible.   Interesting explicit examples in two variables are constructed in [A3]. ---As we shall see, the classical solutions are too few to solve the Dirichlet problem.

\paragraph{Gradient flow} Let $x=x(t) = (x_1(t),x_2(t),\dots,x_n(t))$ be a sufficiently smooth curve. Differentiating $|\nabla u(x(t))|^2$ along the curve, we obtain
$$\dfrac{d}{dt}|\nabla u(x(t))|^2 \,=\, 2\sum_{i,j=1}^{n} \frac{\partial u}{\partial x_j} \frac{\partial^2u}{\partial x_i \partial x_j}\frac{d x_i}{d t}.$$
Select the curve so that it is a solution to the system
$$\frac{d x}{dt} = \nabla u(x(t)).$$
Along this so-called gradient flow, which yields the stream lines\footnote{The stream lines are orthogonal to the level surfaces $u = $ constant.}, we have 
$$\dfrac{d}{dt}|\nabla u(x(t))|^2 \,=\, 2\Delta_{\infty}u(x(t)).$$
This formula has an interesting interpretation for classical solutions of the equation $\Delta _{\infty}u = 0.$  It follows that 
\textsf{the 'speed'} $|\nabla u|$ \textsf{is constant along a stream line} (but different stream lines may have different constants). Recall  from elementary calculus that the level surfaces of a function $u$ have the normal $\nabla u$ and that the growth is at its fastest in the direction of the gradient.

\paragraph{Variational Solutions} Given a Lipschitz continuous function $g$ defined on the boundary of a bounded domain $\Omega$ we first minimize the variational integral 
$$I(u) = \int_{\Omega}\!|\nabla u(x)|^p\,dx$$
among all functions in $C(\overline{\Omega}) \cap C^1(\Omega)$ with boundary values $g.$ There exists a unique minimizer, say $u_p.$ It is a solution of the $p$-Laplace Equation (\ref{pl}). A \emph{variational solution} of the $\infty$-Laplace equation is obtained as a uniform limit
$$u_{\infty} = \lim_{p_j \to \infty}u_{p_j} $$
via suitable subsequences.  The function will be  Lipschitz continuous. Since
$$ \lim_{p \to \infty}\left\{\int_{D}\!|\nabla u(x)|^p\,dx\right\}^{\frac{1}{p}} = \underset{x \in D}{\esssup}\, |\nabla u(x)| =\|\nabla u\|_{L^{\infty}(D)},$$
one may guess that $u_{\infty}$ minimizes the  ''integral'' $\|\nabla u\|_{L^{\infty}}.$ Indeed, for each subdomain $D$ it holds that
$$\|\nabla u_{\infty}\|_{L^{\infty}(D)} \leq \|\nabla v\|_{L^{\infty}(D)},$$
whenever $v$ is Lipschitz continuous in $\overline{D}$ and $v = u_{\infty}$ on the boundary $\partial D.$ Taking into account that 
$$\sup_{D\times D}\dfrac{|u(x)-u(y)|}{|x-y|} = \|\nabla u\|_{L^{\infty}(D)},$$
we interpret it as having obtained an \textsf{''optimal'' Lipschitz extension} of the boundary values $g$ to the domain $\Omega.$ --As we will see, all solutions are variational.

\paragraph{Taylor Expansion and Mean Values} Substituting $h = \pm \varepsilon \frac{\nabla u(x)}{|\nabla u(x)|}$ in the  Taylor expansion
$$u(x+h) = u(x) + \langle\nabla u(x),h\rangle + \frac{1}{2}\langle \mathrm{D}^2u(x)h,h\rangle + \cdots$$
we arrive at the formula
$$\dfrac{u\!\left(x+\varepsilon\frac{\nabla u(x)}{|\nabla u(x)|}\right) + u\!\left(x-\varepsilon\frac{\nabla u(x)}{|\nabla u(x)|}\right)}{2} - u(x) = \frac{1}{2}\varepsilon^2\dfrac{\Delta_{\infty}u(x)}{|\nabla u(x)|^2} + O(\varepsilon^4),$$
which  can be used to develop a numerical scheme\footnote{A function $u \in C^2(\Omega)$  obtained by the so-called Casas-Torres image interpolation algorithm satisfies
$$u(x) = \dfrac{u\!\left(x+h\frac{\nabla u(x)}{|\nabla u(x)|}\right) + u\!\left(x-h\frac{\nabla u(x)}{|\nabla u(x)|}\right)}{2} + o(h^2).$$ See \textsc{J. R. Casas} and \textsc{L. Torres,} \emph{Strong edge features for image coding}, in 'Mathematical Morphology and its Application to Image and Signal Processing' edited by R. Schaefer, P. Maragos, and M. Butt, Kluwer Academic Press, Atlanta 1996, pp. 443--450.}. --Since a  function grows fastest in the direction of its gradient, we see that the ''\emph{Mean Value}'' Formula
$$\dfrac{\underset{\overline{B(x,\varepsilon)}}{\max}\{u(x)\} + \underset{\overline{B(x,\varepsilon)}}{\min}\{u(x)\}}{2}\, -\, u(x)\, =\, \dfrac{1}{2}\varepsilon^2\dfrac{\Delta_{\infty}u(x)}{|\nabla u(x)|^2} + o(\varepsilon^2)$$
holds,  all this provided that $\nabla u(x) \not = 0.$ 
This suggests that $\infty-$harmonic functions can be defined through the asymptotic expansion
$$u(x)= \dfrac{\underset{\overline{B(x,\varepsilon)}}{\max}\{u(x)\} + \underset{\overline{B(x,\varepsilon)}}{\min}\{u(x)\}}{2} + o(\varepsilon^2)\qquad \text{as}\qquad \varepsilon \longrightarrow 0.$$
Properly interpreted this is, indeed, the case. See Section \ref{secasy}.

\paragraph{Comparison with Cones}  An extremely simple property characterizes the $\infty$-harmonic functions: a Comparison Principle with respect to the ''cone'' functions 
$$C(x) = a + b|x-x_0|,$$
which play the role of fundamental solutions (they are $\infty$-harmonic except at their apex $x_0$). The requirement is that, for each subdomain $B$, the inequalities
$$\min_{\xi \in \partial B}\left\{u(\xi)-b|\xi-x_0|\right\} \leq u(x)-b|x-x_0| \leq \max_{\xi \in \partial B}\left\{u(\xi)-b|\xi-x_0|\right\}$$
hold when $x\in B,\, x_0 \in \Rn \setminus B.$
 It is surprising that so small a family of functions suffices. The corresponding theorem for ordinary harmonic functions is more involved (!), see [CZ].

\paragraph{Image Processing} The $\infty$-Laplace operator appears in connexion with image processing, where the $\infty$-Laplace Equation was rediscovered, see [CMS]. It  has several advantages over the ordinary Laplace Equation, which it replaces in some applications. The most useful property is that boundary values can be prescribed even at isolated boundary points (this is out of the question for harmonic functions!). The  boundary can be arbitrarily irregular.  Moreover there is ''no flow along the contours'' due to the fact that the $\infty$-Laplacian is the second derivative in the direction of the stream lines:

$$\dfrac{\Delta_{\infty}u}{|\nabla u|^2} = \frac{\partial^2u}{\partial \nu ^2},\qquad \nu = \dfrac{\nabla u}{|\nabla u|}.$$
For all this we refer to [CMS].

\paragraph{Tug-of-War}  It is well known that a probabilistic approach to the Laplace Equation $\Delta u = 0$ is provided by the celebrated Brownian Motion, also known as the Wiener Process. Quite unexpectedly, the $\infty$-Laplace Equation $\Delta_{\infty}u = 0$ can be given a probabilistic interpretation in terms of \emph{game theory}. A game called Tug-of-War replaces the Brownian motion when it comes to the Dirichlet Boundary value problem. While all directions are equally probable for the Brownian motion, the (mathematical) Tug-of-War
favours only two (generic) directions, viz. $\pm \nabla u$. See Section \ref{sectug}.

\paragraph{Viscosity Solutions}  Last but not least, we mention that if one wants to study the $\infty$-Laplace Equation within the framework of Partial Differential Equations the concept of solutions is delicate. There is not enough of classical solutions and the uniqueness question is difficult if one restricts the attention to variational solutions. The modern theory of \textsf{viscosity solutions} is the right framework. This is the way that the author of these lines prefers. ---Although  an approach, based merely on comparison with cones is possible so that the differential equation  disappears (!), I find that method rather limited in scope. 

\bigskip

There are many variants and extensions of the Infinity-Laplace Operator. A short list is the following.
\begin{itemize}
\item The game theoretic or normalized operator is
$$\Delta_{\infty}^{\diamondsuit}u\, = \,|\nabla u|^{-2}\Delta_{\infty}u.$$
\item With variable coefficients $$
\sum_{i,j=1}^{n}a^{ij}(x)\frac{\partial u}{\partial x_i}\frac{\partial u}{\partial x_j}\frac{\partial^2u}{\partial x_i \partial x_j}\,=\,0.$$
\item The Euler-Lagrange equation $$
\sum_{i,j,k,l}^n\Bigr(2g^{ik}(x)\frac{\partial u}{\partial x_k}\,\frac{\partial^2u}{\partial x_i \partial x_j}\,g^{jl}(x)\frac{\partial u}{\partial x_l} + \frac{\partial g^{ij}(x) }{\partial x_k}\,\frac{\partial u}{\partial x_i}\frac{\partial u}{\partial x_j}\,g^{kl}(x)\frac{\partial u}{\partial x_l}\Bigl)\,=\,0.
$$
for the problem
$$\min_{u}\max_{x} \Bigl|\sum_{i,j=1}^ng^{ij}(x)\frac{\partial u}{\partial x_i}\frac{\partial u}{\partial x_j}\Bigr|.$$
\item The eigenvalue problem
$$\max\Bigl\{\Lambda - \frac{|\nabla u|}{u},\,\sum_{i,j=1}^{n}\frac{\partial u}{\partial x_i}\frac{\partial u}{\partial x_j}\frac{\partial^2u}{\partial x_i \partial x_j}\Bigr\}\, =\, 0.$$
\item The Euler-Lagrange equation
$$\Delta_{\infty}u\,+\, |\nabla u|^2\ln|\nabla u|\,\langle\nabla u,\nabla \log p(x)\rangle \, =\, 0$$
for the problem 
$$\min_{u}\max_{x} |\nabla u(x)|^{p(x)}.$$
\item Time dependent equations  like
$$\frac{\partial u}{\partial t} = \Delta_{\infty}u,\quad u^2\frac{\partial u}{\partial t} = \Delta_{\infty}u,\quad |\nabla u|^2\frac{\partial u}{\partial t}\,=\,\Delta_{\infty}u.  $$
\item Euler-Lagrange equations for $$\min_{u}\max_{x}\bigl\{F(x,u,\nabla u)\bigr\}.$$
\item Discrete versions, so-called $\infty$-harmonious (sic!) functions.
\item The minimax-problem $$\min_{\mathbf{u}}\max_{x} |\mathbb{D} \mathbf{u}(x)|\quad \text{where}\quad
\mathbf{u} = (u^1,u^2,\dots,u^N)$$ yields formally the system
$$\sum_{\mu=1}^N\sum_{i,j=1}^n\frac{\partial u^{\alpha}}{\partial x_i}\frac{\partial u^{\mu}}{\partial x_j}\frac{\partial^2u^{\mu}}{\partial x_i \partial x_j} \, = 0,\quad \alpha = 1,2,\dots, N,$$
for which viscosity solutions do not seem to work.
\end{itemize}
Interesting as they are, we leave these problems aside\footnote{A serious omission is the differentiability proof of the gradient. I have not succeded in simplifying the proofs in [S] and [ES], and so it is best to directly refer to the original publications.} and stick to the more pregnant original formulation. ---''Less is more''.

\section{Preliminaries}
\label{secpre}
 \texttt{Abundant dulcibus vitiis}

\bigskip
 
 \noindent I immediately mention that always\footnote{The $p$-harmonic operator is not pointwise defined for $p <2.$} $p \geq 2$ in these notes and often $p > n =$ the dimension of the space\footnote{All functions in the Sobolev space $W^{1,p}$ are continuous when $p>n.$}. (At a first reading, one had better keep $p$ large.) We aim at $p = \infty.$

We begin with a special case of Ascoli's theorem. 

\begin{thm}[Ascoli] Let the sequence of functions $f_k : \Omega \longrightarrow \R$ be
\begin{description}
\item equibounded: $$\sup_{\Omega}|f_{k}(x)|\, \leq M\, < \infty \quad \text{ when}\quad k=1,2,\dots,$$
\item equicontinuous: $$|f_k(x)-f_k(y)| \,\leq\, C|x-y|^{\alpha}\quad \text {when} \quad k=1,2,\dots.$$
\end{description}
Then there exists a continuous function $f$ and a subsequence such that\\ $f_{k_j} \longrightarrow f$ locally uniformly in $\Omega$.

 If the domain $\Omega$ is bounded, all functions can be extended continuously to the boundary, and the convergence is uniform in the closure $\overline{\Omega}.$
\end{thm}

\bigskip

{\small
\emph{Proof:} We reproduce a well-known proof. First, we construct a subsequence which converges at the rational points. Let $q_1,q_2,q_3,\dots $ be a numbering of the rational points in $\Omega.$ Since the sequence $f_1(q_1),f_2(q_1),f_3(q_1),...$
is bounded  by our assumption it has a convergent subsequence (Weierstrass' Theorem), say $f_{1j}(q_1),\, j = 1,2,3,\dots.$ Consider the next point $q_2.$ Now the sequence $f_{11}(q_2),\\f_{12}(q_2),f_{13}(q_2),\dots $ is bounded and so we can  extract a convergent subsequence, say $f_{21}(q_2),f_{22}(q_2),f_{23}(q_2),\dots $. Continuing like this, extracting subsequences of subsequences, we have the scheme
\begin{align*}
&f_{11},f_{12},f_{13},f_{14}\dots  \quad \text{converges at}\quad q_1\\
&f_{21},f_{22},f_{23},f_{24}\dots  \quad \text{converges at}\quad q_1,q_2\\
&f_{31},f_{32},f_{33},f_{34}\dots  \quad \text{converges at}\quad q_1,q_2,q_3\\
&f_{41},f_{42},f_{43},f_{44}\dots  \quad \text{converges at} \quad q_1,q_2,q_3,q_4\\
&\cdots\cdots\cdots\cdots\cdots \cdots\cdots\cdots\cdots\cdots \cdots\cdots\cdots \cdots\cdots.
\end{align*}
We see that the diagonal sequence $f_{11},f_{22},f_{33},f_{44},\dots$ converges at every rational point. To simplify the notation, write $f_{k_j} = f_{jj}.$

 We claim that the constructed diagonal sequence converges at each point in $\Omega,$ be it rational or not. To this end let $x\in \Omega$ be an arbitrary point and let $q$ be a rational point very near it.
Then
\begin{align*}
|f_{k_j}(x)-f_{k_i}(x)| &\leq |f_{k_j}(x)-f_{k_j}(q)| + |f_{k_j}(q)- f_{k_i}(q)| + |f_{k_i}(q)-f_{k_i}(x)|\\&\leq 2C|x-q|^{\alpha} + |f_{k_j}(q)- f_{k_i}(q)|.
\end{align*}
Given $\varepsilon > 0,$ we fix $q$ so close to $x$ that $2C|x-q|^{\alpha}< \tfrac{\varepsilon}{2},$ which is possible since the rational points are dense. By the convergence at the rational points, we infer that 
$$|f_{k_j}(x)-f_{k_i}(x)| < \frac{\varepsilon}{2} + \frac{\varepsilon}{2} = \varepsilon,$$
when the indices $i$ and $j$ are large enough. By Cauchy's general convergence criterion, the sequence converges at the point $x.$ We have established the existence of the pointwise limit function
$$f(x) = \lim_{j \to \infty} f_{k_j}(x).$$

Next, we show that the convergence is (locally) uniform.
Suppose that $\overline{\Omega}$ is compact. Cover it by balls $B(x,r)$
with diameter $2r = \varepsilon^{\frac{1}{\alpha}}.$
A finite number of these balls covers  $\overline{\Omega}$:
$$\overline{\Omega}\, \subset\, \bigcup_{m=1}^{N}B(x_m,r).$$ Choose a rational point from each ball, say $q'_m \in B(x_m,r).$ Since  only a finite number of these points are involved, we can fix an index $N_{\varepsilon}$ such that
$$\max_{m}|f_{k_j}(q'_m)- f_{k_i}(q'_m)| < \varepsilon \qquad \text{when} \quad i,j > N_{\varepsilon} $$ 
Let $x\in \overline{\Omega}$ be arbitrary. It must belong to some ball, say $B(x_m,r).$
Again we can write
\begin{align*}
|f_{k_j}(x)-f_{k_i}(x)|&\leq 2C|x-q'_m|^{\alpha} +|f_{k_j}(q'_m)- f_{k_i}(q'_m)| \\&\leq
2C(2r)^{\alpha} +|f_{k_j}(q'_m)- f_{k_i}(q'_m)|\\ &\leq 2C\varepsilon +\varepsilon \qquad \text{when} \quad i,j > N_{\varepsilon}. 
\end{align*}
The index $N_{\varepsilon}$ is independent of how the point $x$ was chosen. This shows that the convergence is \emph{uniform} in  $\overline{\Omega}.$ ---If the domain $\Omega$ is unbounded, we notice that the above proof is valid for every fixed bounded subdomain, in which case the concergence is locally uniform. \qquad $\Box$}

\bigskip

Next, we consider Lipschitz continuous functions. A function $f : \Omega \longrightarrow \R$ is Lipschitz contiunuous if for some constant $L$,
$$|f(x)-f(y)| \leq L|x-y|\quad\text{when}\quad x,y \in \Omega.$$

\begin{thm}[Rademacher] A Lipschitz continuous function $f$ is totally differentiable a.\,e. in its domain: the expansion
$$f(y) = f(x) + \langle \nabla f(x),y-x\rangle + o(|y-x|)\quad\text{as}\quad y \to x$$
holds at almost every point $x \in \Omega.$
\end{thm}

\bigskip

It is useful to know that convex functions are locally Lipschitz continuous. As we shall see in Section \ref{secfro}, a convex function has, indeed, even second derivatives  a.\,e. in the way they should appear in Taylor's expansion.

\paragraph{Sobolev Spaces} We denote by $W^{1,p}(\Omega)$ the Sobolev space consisting of functions $u$ that together with their first distributional derivatives
$$\nabla u = \Bigl(\frac{\partial u}{\partial x_1},\frac{\partial u}{\partial x_2},\dots,\frac{\partial u}{\partial x_n}\Bigr)$$
belong to the space $L^p(\Omega).$ Equipped with the norm
$$ \|u\|_{W^{1,p}(\Omega)}\, = \, \|u\|_{L^p(\Omega)} + \|\nabla u\|_{L^p(\Omega)}$$
it is a Banach space. In particular, the space $W^{1,\infty}(\Omega)$ consists of all Lipschitz continuous functions defined in $\Omega.$ The closure of $C^{\infty}_0(\Omega)$ with respect to the Sobolev norm is denoted by $W^{1,p}_0(\Omega).$  About Sobolev spaces we refer the reader to Chapter 7 of the book [GT].

 Before proceeding, let us take advantage of the fact that only very large values of the exponent $p$ are needed here. If $p > n =$ the number of coordinates in $\Rn,$ the Sobolev space contains only continuous functions and the boundary values are taken in the classical sense. \textsf{All domains\footnote{As always, a domain is an open connected set.} $\Omega$ are regular for the Dirichlet problem, when $\mathbf{p > n}.$}

\begin{lemma}
\label{Mor}
 Let $p > n$ and suppose that $\Omega$ is an arbitrary bounded domain in $\Rn.$ If $v \in W^{1,p}_0(\Omega),$ then 
\begin{equation}
|v(x)-v(y)|\,\leq\, \frac{2pn}{p-n}\,|x-y|^{1-\frac{n}{p}}\|\nabla v\|_{L^p(\Omega)}
\end{equation}
for a.\,e. $x,y \in \Omega.$ One can redefine $v$ in a set of measure zero and extend it to the boundary so that $v \in C^{1-\frac{n}{p}}(\overline{\Omega})$ and $v|_{\partial \Omega} = 0.$
\end{lemma}

This is a variant of Morrey's inequality. It is important that the constant 
 remains bounded for large $p.$ If we do not require zero boundary values, the inequality still holds for many domains. For example, if $\Omega$ is a cube $Q$, the inequality holds for $v\in W^{1,p}(Q).$

\bigskip

{\small \paragraph{On the Constant} Since the behaviour of the constant is decisive, as $p \to \infty,$ I indicate how to obtain it for a smooth function $v \in C^1(Q) \cap W^{1,p}(Q).$ Integrating
\begin{align*}
v(x)-v(y) =& \int_0^1\!\frac{d}{dx}\,v\bigl(x+t(y-x)\bigr)\,dt\\
 =&  \int_0^1\!\bigl\langle y-x,\nabla v\bigl(x+t(y-x)\bigr)\bigr\rangle\,dt
\end{align*}
with respect to $y$ over $Q$, we see that\footnote{The notation $$f_Q = \vint_Q\!f\,dx = \frac{\int_Q\!f\,dx}{\int_Q\!dx}$$ is used for the average of a function.}
\begin{align*}
|v(x)-v_Q| &=\left|\vint_Q\!\int_0^1\!\bigl\langle y-x,\nabla v\bigl(x+t(y-x)\bigr)\bigr\rangle\,dt\,dy\right| \\&\leq \mathrm{diam}(Q)\vint_Q\!\int_0^1\!|\nabla v\bigl(x+t(y-x)\bigr)|\,dt\,dy\\ &\leq
\mathrm{diam}(Q)\int_0^1\!\Bigl(\vint_Q|\nabla  v\bigl(x+t(y-x)\bigr)|^p\,dy\Bigr)^{\frac{1}{p}}dt.
\end{align*}
The change of coordinates
$$\xi = x + t(y-x),\quad d\xi = t^ndy$$
in the inner integral yields 
$$\vint_Q|\nabla  v\bigl(x+t(y-x)\bigr)|^p\,dy = \frac{1}{t^n}\vint_{Q_t}|\nabla  v(\xi)|^p\,d\xi \leq  \frac{1}{t^n}\vint_{Q}|\nabla  v(\xi)|^p\,d\xi,$$
since the intermediate domain of integration $Q_t \subset Q.$
Therefore
$$|v(x)-v_Q| \leq \frac{\mathrm{diam}(Q)}{|Q|^\frac{1}{p}}\int_0^1\!t^{-\frac{n}{p}}\|\nabla v\|_{L^p(Q)}\,dt =
\dfrac{1}{1-\frac{n}{p}}\frac{\mathrm{diam}(Q)}{|Q|^{\frac{1}{p}}}\|\nabla v\|_{L^p(Q)}.$$
(It was needed that $p > n.$) The triangle inequality yields
$$|v(x)-v(y)| \leq |v(x)-v_Q| + |v(y) -v_Q| \leq \frac{2p}{p-n}\frac{\mathrm{diam}(Q)}{|Q|^{\frac{1}{p}}}\|\nabla v\|_{L^p(Q)}.$$
To conclude, we can always choose an auxiliary cube $Q' \subset Q$ so that $|x-y| \leq \mathrm{diam}(Q').$ ---In the general case, when $v$ no longer has continuous first derivatives, one  can use approximations with convolutions and conclude the proof with the aid of Ascoli's theorem.}

\bigskip

 I repeat  that always\footnote{The $p$-harmonic operator is not pointwise defined for $p <2.$} $\mathbf{p \geq 2}$ in these notes and often $\mathbf{p > n} =$ the dimension of the space\footnote{All functions in the Sobolev space $W^{1,p}$ are continuous when $p>n.$}. 

\paragraph{The $p$-Laplace Equation for finite $p$} We need some standard facts about the $p$-Laplace equation and its solutions. Let us consider the existence of a solution to the Dirichlet boundary value problem. Minimizing the variational integral
\begin{equation}
\label{energy}
I(v) = \int_{\Omega}\!|\nabla v|^p\,dx
\end{equation}
among all functions with the same given boundary values, we are led to the condition that the first variation must vanish:
\begin{equation}
\label{fvariation}
\int_{\Omega}\!\langle |\nabla u|^{p-2}\nabla u,\nabla \eta\rangle\,dx\,=\,0\quad\text{when}\quad \eta \in C^{\infty}_{0}(\Omega),
\end{equation}
where $u$ is minimizing. Under suitable assumptions this is equivalent to
$$\int_{\Omega}\!\eta\,\nabla\!\cdot\!\bigl(|\nabla u|^{p-2}\nabla u\bigr)\,dx\,=\,0,$$
Since this  must hold for all  test functions $\eta$ we  have 
$$\Delta_{p}u\,\equiv\nabla\!\cdot\!\bigl(|\nabla u|^{p-2}\nabla u\bigr)\,=\,0.$$
In other words, the $p$-Laplace Equation is the \emph{Euler-Lagrange Equation} for the above variational integral. A more precise statement is:

\bigskip

\begin{thm}
\label{Lagrange}
Take $p > n$ and consider an arbitrary bounded domain  $\Omega$ in $\Rn.$ Suppose that $g\in C(\overline{\Omega})\cap W^{1,p}(\Omega)$ is given. Then there exists a unique function $u\in C(\overline{\Omega})\cap W^{1,p}(\Omega)$ with boundary values $g$ which minimizes the variational integral 
$$I(v) = \int_{\Omega}|\nabla v|^p\,dx$$
 among all similar functions.

 The minimizer is a weak solution to the $p$-Laplace Equation, i.e..
$$\int_{\Omega}\!\langle |\nabla u|^{p-2}\nabla u,\,\nabla \eta\rangle\,dx\,=\,0\quad\text{when}\quad \eta \in C^{\infty}_{0}(\Omega).$$ On the other hand, a weak solution in $C(\overline{\Omega})\cap W^{1,p}(\Omega)$ is always a minimizer (among functions with its own boundary values).
\end{thm}

\bigskip

\emph{Proof:} The uniqueness of the minimizer follows easily from 
$$\left|\dfrac{\nabla u_1 + \nabla u _2}{2}\right|^p\, < \dfrac{|\nabla u_1|^p+|\nabla u_2|^p}{2}\quad \text{when}\quad \nabla u_1\not = \nabla u_2,$$
upon integration. To wit, if there were two minimizers $u_1$ and $u_2,$ then $\tfrac{u_1+u_2}{2}$ would be admissible and
$$I(u_1)\leq I\Bigl(\frac{u_1+u_2}{2}\Bigr) < \dfrac{I(u_1)+I(u_2)}{2} = I(u_1),$$
unless $\nabla u_1 = \nabla u_2$ almost everywhere. To avoid the contradiction, we must have $u_1=u_2.$

The Euler-Lagrange Equation can be derived from the minimizing property
$$I(u) \leq I(u+\varepsilon \eta).$$
(The function $v(x) = u(x) + \varepsilon \eta(x)$ is admissible.)  We must have
$$\dfrac{d}{d\varepsilon}I(u+\varepsilon \eta) \,=\,0\quad\text{when}\quad\varepsilon =0$$
by the infinitesimal calculus.
This shows that the first variation vanishes, i.e, equation (\ref{fvariation}) holds.

To show that the minimizer exists, we use the Direct Method in the calculus of variations, due to Lebesgue, see the book [G]. Let
$$I_0 \,=\, \inf_{v} \int_{\Omega}\!|\nabla v|^p\,dx$$
where the infimum is taken over the class of admissible functions. Now $0\leq I_0\leq I(g) < \infty.$ Consider a so-called minimizing sequence of admissible functions $v_j$:
$$\lim_{j\to\infty}I(v_j)\,=\,I_0.$$
We may assume that $I(v_j) < I_0+1$ for $j=1,2,\dots.$ We may also assume that $$\min g \leq v_j(x) \leq \max g\quad\text{ in}\quad \Omega,$$ since we may cut the functions at the constant  heights $\min g$ and $\max g.$ (The procedure decreases the integral!) We see that the Sobolev norms $\|v_j\|_{W^{1,p}}(\Omega)$ are uniformly bounded\footnote{The conventional way is not to cut the functions, but to use the Sobolev inequality $$\|v_j-g\|_{L^p(\Omega)}\leq C\|\nabla(v_j-g)\|_{L^p(\Omega)}$$ to uniformly bound the norms
\begin{align*}
\|v_j\|_{L^p(\Omega)}&\leq \|v_j-g\|_{L^p(\Omega)}+\|g\|_{L^p(\Omega)}\leq C\|\nabla(v_j-g)\|_{L^p(\Omega)}+\|g\|_{L^p(\Omega)}\\&\leq C\bigl[\|\nabla v_j\|_{L^p(\Omega)} +\|\nabla g\|_{L^p(\Omega)}\bigr]+\|g\|_{L^p(\Omega)}\\&\leq C\bigl[(I_0+1)^p +\|\nabla g\|_{L^p(\Omega)}\bigr] +\|g\|_{L^p(\Omega)} \equiv M < \infty,
\end{align*}
 when $j=1,2,3,\dots.$}   . By weak compactness, there exists a function $u \in W^{1,p}(\Omega)$ and a subsequence such that
$$\nabla v_{j_k} \rightharpoonup \nabla u\quad\text{weakly in}\quad L^{p}(\Omega).$$
Since $p > n$ we know from Lemma \ref{Mor} that $u\in C(\overline{\Omega})$ and we conclude that $u=g$ on $\partial \Omega.$ We got the continuity for free! By the weak lower semicontinuity of convex integrals
$$I(u) \leq \liminf_{k\to \infty}I(v_{j_k})\,=\,I_0.$$
Since $u$ is admissible also $I(u) \geq I_0.$ Therefore $u$ is a minimizer and the existence is established.

It remains to show that the weak solutions of the Euler-Lagrange Equation are minimizers\footnote{There are variational integrals for which this is not the case.}. By integrating the inequality\footnote{Since $|w|^p$ is convex, the inequality
$$|b|^p \geq|a|^p + p\langle|a|^{p-2}a,b-a\rangle$$
holds for vectors.}
$$|\nabla(u+\eta)|^p \geq |\nabla u|^p + p\langle|\nabla u|^{p-2}\nabla u,\nabla \eta\rangle,$$
we obtain
$$\int_{\Omega}\!|\nabla(u+\eta)|^p\,dx \geq\int_{\Omega}\!|\nabla u|^p\,dx + 0 =   \int_{\Omega}\!|\nabla u|^p\,dx .$$
Therefore $u$ is a minimizer.
 \qquad $\Box$

\bigskip
\emph{Remark:} If the given boundary values $g$ are merely continuous ($g\in C(\partial \Omega$)  but perhaps $g\not\in W^{1,p}(\Omega)$), then  there exists a unique $p$-harmonic function $u \in C(\overline{\Omega})$ with boundary values $g.$
However it may so happen that $\int_{\Omega}\!|\nabla u|^p\,dx = \infty.$ ---Hadamard gave a counter example for the ordinaruy Dirichlet integral ($p=2$). 

\section{Variational Solutions}
\label{secvar}

The passage to the limit in the $p$-Laplace Equation $\Delta_pu = 0$ as $p \to \infty$ is best done in connexion with the problem of extending Lipschitz boundary values to the whole domain. To this end, let $g:\,\partial\Omega \to 
\R$ be Lipschitz continuous:
$$|g(\xi_1) - g(\xi_2)|\,\leq\,L|\xi_1-\xi_2|,\qquad \xi_1,\xi_2 \in \partial \Omega.$$
If  $u:\,\overline{\Omega} \to 
\R$ is a  Lipschitz continuous function with constant $L$ and boundary values $u|_{\partial \Omega} = g,$ then some elementary manipulations show that necessarily
$$\max_{\xi \in \partial \Omega}\{g(\xi)-L|x-\xi|\}\, \leq u(x) \leq\, \min_{\xi \in \partial \Omega}\{g(\xi)+L|x-\xi|\}.$$
As functions of $x$, the two bounds are Lipschitz continuous extensions of $g$ and have Lipschitz constant $L$. They were found by McShane and Whitney. In several variables they do usually not coincide, uniqueness fails. We aim at finding the best
extension.

By Rademacher's Theorem a Lipschitz function is a.\,e. differentiable. Extending the given $g$, first in some convenient way (for example as the above majorant), we now have
$$\|\nabla g\|_{L^{\infty}(\Omega)}\,\leq\,L,\quad g \in C(\overline{\Omega})\cap W^{1,\infty}(\Omega),$$
 where the same notation is used for the extended function.
We shall minimize the variational integral
$$I_p(v) = \int_{\Omega}\!|\nabla v|^p\,dx$$
among all functions in $W^{1,p}(\Omega)$ with boundary values $v-g \in W^{1,p}_0(\Omega).$
 Then we shall let $p \to \infty.$

 Before proceeding, let us take advantage of the fact that only very large values of the exponent $p$ are needed here. If $p > n =$ the number of coordinates in $\Rn,$ the Sobolev space contains only continuous functions and the boundary values are taken in the classical sense. \textsf{All domains\footnote{As always, a domain is an open connected set.} $\Omega$ are regular for the Dirichlet problem, when $p > n.$}

 By \emph{the direct method in the Calculus of Variations} a unique minimizer $u_p$ exists.  See Theorem \ref{Lagrange}. In particular $u_p \in C(\overline{\Omega})\cap W^{1,p}(\Omega)$ and $u_p = g$ on the boundary $\partial \Omega.$ By the minimization property
$$\|\nabla u_p\|_{L^p(\Omega)} \leq\|\nabla g\|_{L^p(\Omega)}  \leq L|\Omega|^{\frac{1}{p}}.$$

 We aim at using Ascoli's Theorem.
Morrey's inequality in Lemma \ref{Mor}  is applicable to the function $v_p = u_p - g.$ Thus we have 
\begin{align*}
|u_p(x)&-u_p(y)|\\&\leq |g(x)-g(y)| + |v_p(x)-v_p(y)|\\
&\leq L|x-y| + \frac{2pn}{p-n}\,|x-y|^{1-\frac{n}{p}}\|\nabla v_p\|_{L^p(\Omega)}\\
&\leq L|x-y| + \frac{2pn}{p-n}\,|x-y|^{1-\frac{n}{p}}(\|\nabla u_p\|_{L^p(\Omega)}+\|\nabla g\|_{L^p(\Omega)}) \\
&\leq  L|x-y| + LC_{\Omega}|x-y|^{1-\frac{n}{p}},
\end{align*}
when $p > n+1,$ say. Also the bound 
$$\|u_p\|_{L^{\infty}(\Omega)} \leq  \|g\|_{L^p(\Omega)} + LC'_{\Omega}$$
follows. Thus the family $\{u_p\},\,\, p > n+1,$ is \emph{equicontinuous} and \emph{ equibounded.} Therefore Ascoli's theorem applies and we can extract a subsequence $p_j \to \infty$ such that
$$u_{p_{j}} \longrightarrow u_{\infty} \quad \text{uniformly in}\quad \overline{\Omega},$$
where $u_{\infty} \in  C(\overline{\Omega})$ is some function with boundary values $u_{\infty} |_{\partial \Omega} = g.$ Moreover,
$$|u_{\infty}(x) - u_{\infty}(y)|\, \leq\, C\,|x-y|,$$
where $C$ depends on $L$ and $\Omega.$ Observe that the limit function may depend on how the sequence $p_j$ was extracted!

We can do better than this. When $p_j > s,$ H\"{o}lder's inequality yields\footnote{The convenient notation
$$\vint_{A}\!f\,dx = \frac{\int_{A}\!f\,dx}{\int_{A}\,dx}$$ 
is used for the average of a function. Then $\{\vint_{A}\!|f|^p\,dx\}^{\frac{1}{p}}$ is increasing with $p.$}
$$\left\{\vint_{\Omega}\!|\nabla u_{p_j}|^s\,dx\right\}^{\frac{1}{s}} \leq \left\{\vint_{\Omega}\!|\nabla u_{p_j}|^{p_j}\,dx\right\}^{\frac{1}{p_j}} \leq \left\{\vint_{\Omega}\!|\nabla g|^{p_j}\,dx\right\}^{\frac{1}{p_j}} \leq L.$$
Thus
$$\nabla u_{p_{j}}\, \rightharpoonup\, \nabla u_{\infty} \quad \text{weakly in}\quad L^s(\Omega)$$ for some subsequence of $p_j.$ By weak lower semicontinuity
$$\int_{\Omega}\!|\nabla u_{\infty}|^s\,dx \leq \liminf_{j \to \infty}\int_{\Omega}\!|\nabla u_{p_j}|^s\,dx.$$
We remark that, via a diagonalization procedure\footnote{It is helpful that weak convergence in $L^p$ implies weak convergence in $L^s,\, s< p.$}, it is possible to extract a single subsequence so that
$$\nabla u_{p_{j}}\, \rightharpoonup\, \nabla u_{\infty} \quad \text{weakly in each}\quad L^s(\Omega)\quad \text{simultaneously},$$ 
where $n+1 < s < \infty$.
It follows that
$$\left\{\vint_{\Omega}\!|\nabla u_{\infty}|^s\,dx\right\}^{\frac{1}{s}} \leq L.$$
But $s$ was arbitrarily large. Letting $s \to \infty$ we conclude that
\begin{equation}
\|\nabla u_{\infty}\|_{L^{\infty}(\Omega)} \leq L.
\end{equation}
Thus the Lipschitz constant of $u_{\infty}$ is $L.$\footnote{This is the case also for non-convex domains, because $|g(\xi_1)-g(\xi_2)| \leq L|\xi_1 - \xi_2|$ even if the segment joining the boundary points $\xi_1$ and $\xi_2$ is outside $\Omega.$}

 A difficult problem remains, namely the uniqueness of $u_{\infty}.$ Does the limit depend on the particular subsequence? Actually, it is unique, but no proof using only variational methods is known to us. Other tools are required. We postpone the issue and record an immediate result, according to which one cannot even \emph{locally} improve the Lipschitz constant. It is the best one!

\begin{thm}[Existence]
\label{var}
Given $g \in C(\overline{\Omega}) \cap W^{1,\infty}(\Omega),$ there exists a function $u_{\infty} \in C(\overline{\Omega}) \cap W^{1,\infty}(\Omega)$ with boundary values $u_{\infty} = g$ on $\partial \Omega$ having the following minimizing property in each subdomain $D \subset \Omega:$

If $v \in C(\overline{D}) \cap W^{1,\infty}(D),$ and $v = u_{\infty}$ on $\partial D,$ then
$$\|\nabla u_{\infty}\|_{L^{\infty}(D)} \leq \|\nabla v\|_{L^{\infty}(D)}.$$
This $u_{\infty}$ can be obtained as the uniform limit $\lim u_{p_j}$ in $\Omega,$ \,$u_{p_j}$ denoting the solution of the $p_j$-Laplace equation with boundary values $g.$
\end{thm}

\emph{Proof:} It remains to show that the constructed
$$u_{\infty} = \lim_{j \to \infty} u_{p_j}$$
has the minimizing property in $D.$ Since the uniqueness is not yet available to us, we have to be careful. Let $v$ be the given function in $\overline{D}$ and let $v_{p_j}$ denote the solution to $\Delta_{p_j}v_{p_j} = 0$ in $D$ with boundary values $v_{p_j}|_{\partial D} = u_{\infty}|_{\partial D}.$ Then $v_{p_j}$ has the minimizing property
$$\int_{\Omega}\!|\nabla v_{p_j}|^{p_j}\,dx \leq \int_{\Omega}\!|\nabla v|^{p_j}\,dx.$$

We claim that $v_{p_j} \to u_{\infty}$ uniformly in $\overline{D}.$ To this end, notice that since
$$v_{p_j}- u_{\infty} = (v_{p_j}-u_{p_j})+(u_{p_j}-u_{\infty})$$
it is sufficient to show that $\|v_{p_j}-u_{p_j}\|_{L^{\infty}(D)} \longrightarrow 0.$ Now
\begin{align*}
&\max_{\overline{D}}(v_{p_j}-u_{p_j}) = \max_{\partial D}(v_{p_j}-u_{p_j})\\
= &\max_{\partial D}(u_{\infty}-u_{p_j}) \leq \|u_{\infty}-u_{p_j}\|_{L^{\infty}(D)}\to0
\end{align*}
and the same for $u_{p_j}-v_{p_j}.$ This proves the claim.  {\small ---We used the fact that the difference of two solutions of the $p$-Laplace Equation attains its maximum at the boundary. To see this, assume that
$$\alpha = \max_{\overline{D}}(v_{p}-u_{p}) >  \max_{\partial D}(v_{p}-u_{p}) = \beta.$$
Define the open set
$$G = \left\{x\in D|\,v_{p}(x)-u_{p}(x) > \tfrac{\alpha+\beta}{2}\right\}.$$
Then $G \subset \subset D$ and $v_p = u_p+\tfrac{\alpha+\beta}{2}$ on $\partial G.$
But both $v_p$ and $u_p+\tfrac{\alpha+\beta}{2}$ are solutions of the $p$-Laplace Equation in $G$ with the same boundary values. By uniqueness they coincide in $G,$ which is a contradiction at the maximum point(s). Thus the maximum was at the boundary, as claimed.}

Let us proceed, knowing that $v_{p_j} \to u_{\infty}$ in $D$. Fix $s > > 1.$ Again
$$\left\{\vint_D\!|\nabla v_{p_j}|^s\,dx\right\}^{\frac{1}{s}} \leq \left\{\vint_D\!|\nabla v_{p_j}|^{p_j}\,dx\right\}^{\frac{1}{p_j}} \leq \left\{\vint_D\!|\nabla v|^{p_j}\,dx\right\}^{\frac{1}{p_j}} \leq \|\nabla v\|_{L^{\infty}(D)}$$
when $p_j > s.$
We can extract a weakly convergent subsequence $\nabla v_{p_j} \rightharpoonup \nabla u_{\infty}$ in $L^s(D)$ and by weak lower semicontinuity
$$\left\{\vint_D\!|\nabla u_{\infty}|^{s}\,dx\right\}^{\frac{1}{s}} \leq \|\nabla v\|_{L^{\infty}(D)}.$$
Let $s \to \infty$ to finish the proof. $\Box$

\bigskip

As an example of a property that can be deduced via the variational procedure we consider Harnack's Inequality.

\begin{thm}[Harnack] If the variational solution $u_{\infty} \geq 0$ in $\Omega$, then
$$\boxed{u_{\infty}(x) \leq u_{\infty}(y)\exp\Bigl(\frac{|x-y|}{R-r}\Bigr)}$$
when $x,y \in B(x_0,r)$ and $\overline{B(x_0,R)} \subset \Omega.$
\end{thm}

\emph{Proof:} Let $u_{\infty} = \lim u_p$ via some subsequence. Then $u_{\infty}+\varepsilon = \lim( u_p+\varepsilon)$ and since constants may be added to solutions, we may assume that $u_p \geq \varepsilon > 0.$ At the end, $\varepsilon$ can be sent to $0.$

Let $\zeta \in C^{\infty}_0(\Omega),\,0\leq \zeta \leq 1$ and use the test function
$$\eta = \frac{\zeta^p}{u_p^{p-1}}$$
in the equation
$$\int_{\Omega}\!\langle|\nabla u_p|^{p-2}\nabla u_p,\nabla \eta \rangle \,dx = 0.$$
From
$$\nabla \eta = p\Bigl(\frac{\zeta}{u_p}\Bigr)^{p-1}\nabla \zeta - (p-1)\Bigl(\frac{\zeta}{u_p}\Bigr)^{p}\nabla u_p$$
we obtain
\begin{gather*}
(p-1)\int_{\Omega}\!\zeta^p|\nabla \log u_p|^p\,dx \leq p\int_{\Omega}\!\Bigl(\frac{\zeta}{u_p}\Bigr)^{p-1}|\nabla u_p|^{p-1}|\nabla \zeta|\,dx\\
\leq p\Bigl(\int_{\Omega}\!|\nabla \zeta|^p\,dx\Bigr)^{\frac{1}{p}} \Bigl(\int_{\Omega}\!\zeta^p|\nabla \log u_p|^p\,dx\Bigr)^{1-\frac{1}{p}}
\end{gather*}
and hence
$$\Bigl(\int_{\Omega}\!\zeta^p|\nabla \log u_p|^p\,dx\Bigr)^{\frac{1}{p}} \leq \frac{p}{p-1}\Bigl(\int_{\Omega}\!|\nabla \zeta|^p\,dx\Bigr)^{\frac{1}{p}}.$$
Letting $p \to \infty$ via the subsequence one obtains\footnote{We have skipped the procedure with the auxiliary intermediate space $L^s$ as $s \to \infty$.}
$$\|\zeta \nabla \log u_{\infty}\|_{L^{\infty}(\Omega)} \leq \|\nabla \zeta\|_{L^{\infty}(\Omega)}.$$
Choose $\zeta$ so that $\zeta = 1$ in $B(x_0,r)$ and $|\nabla \zeta| \leq \tfrac{1}{R-r},$ and $\zeta = 0$ outside $B(x_0,R).$ Then
$$\| \nabla \log u_{\infty}\|_{L^{\infty}(B(x_0,r)}\, \leq\, \frac{1}{R-r}.$$
The result follows from 
$$|\log u_{\infty}(x)-\log u_{\infty}(y)|\, \leq \,\|\nabla \log u_{\infty}\|_{L^{\infty}(B(x_0,r)}|x-y|.\qquad \Box$$

 \bigskip

Notice the important conclusion that, if $u_{\infty} \geq 0$ but $u_{\infty} \not \equiv 0,$ then $u_{\infty} > 0$ in $\Omega.$ The virtue of the proof was the good constant.
We remark that the above proof also works for non-negative variational \emph{super}solutions. ---The variational method is used in [GM] to derive Phragm\'{e}n-Lindel\"{o}f Theorems.

{\small \paragraph{Uniqueness for Classical Solutions without Critical Points} It is instructive to prove that 
$$\min_{\Omega}(v-u) \geq \min_{\partial \Omega}(v-u)$$
for $u,v \in C(\overline{\Omega})\cap C^2(\Omega)$ satisfying $\Delta_{\infty}v \leq 0$ and $\Delta_{\infty}u \geq 0$ pointwise in the bounded domain $\Omega$ under the \textsf{extra assumption} $\nabla v \not = 0$. The absence of critical points simplifies the proof. The proof is indirect starting with the antithesis
$$\min_{\Omega}(v-u)  < \min_{\partial \Omega}(v-u).$$
At the interior minimal point we have $\Delta_{\infty}v \geq \Delta_{\infty}u,$ but unfortunately this is not yet a contradiction. Therefore we consider
$$w = \dfrac{1-e^{-\alpha v}}{\alpha} = v - \dfrac{1}{2}\alpha v^2 + \cdots ,$$
where we take $\alpha > 0$ so small that also
$$\min_{\Omega}(w-u)  < \min_{\partial \Omega}(w-u).$$
At the minimal point we have $\Delta_{\infty}w \geq \Delta_{\infty}u \geq 0.$ A calculation shows that
$$\Delta_{\infty}w = e^{-3\alpha v}\left(\Delta_{\infty}v - \alpha |\nabla v|^4\right) \leq  -\alpha  e^{-3\alpha v}  |\nabla v|^4.$$
The extra assumption yields the strict inequality $\Delta_{\infty}w < 0$ at the minimal point. Now this is a contradiction. In particular, this shows uniqueness for the Dirichlet boundary value problem. ---If the extra assumption is abandoned, the construction is more elaborate. In Jensen's proof in Section \ref{secuni} an additional auxiliary equation is used to avoid critical points. Without second derivatives also a 'doubling of variables' is needed.} 

\paragraph{Inheritance Lost} Many properties of the $p$-harmonic functions are passed over to their limits as $p$ goes to $\infty,$ but some qualities are inevitably lost in the passage. The problem of \emph{unique continuation} is actual. Can two $p$-harmonic functions, defined in the same domain, coincide in an open subdomain without being identical?
For a finite $p$ the Principle of Unique Continuation is known to hold in the plane, in space the problem seems to be open. For the $\infty$-Laplace equation the principle is plainly false. This is exhibited by the example
\begin{equation*}
u(x,y) =
\begin{cases}
1-\sqrt{x^2+y^2},\qquad \text{if}\quad x\leq 0,\,y\geq0\\
1-y, \qquad\qquad\quad \,\,\,\,\text{if}\quad x\geq 0,\,y\geq 0
\end{cases}
\end{equation*}
defined in the half-plane   $y\geq 0.$ In other words, we see a cone and its tangent plane, when we draw the surface $z = u(x,y).$ This is a viscosity solution of the $\infty$-Laplace Equation. Both expressions
$$1-\sqrt{x^2+y^2}\qquad \text{and} \qquad 1-y$$
provide $\infty$-harmonic functions that coincide with $u$ in  some subdomains. ---The  example is easily modified to several dimensions.

The following variant of the problem is unsettled, so far as I know. Can an $\infty$-harmonic function vanish at each point in an open  subdomain without being identically zero? The above example does not address this question.

Another phenomenon is the lost \emph{real analyticity}. For a finite $p$ it is known that
the $p$-harmonic functions are of class \,$C^{1,\alpha}_{loc}$,\, i.e., their gradients are locally H\"{o}lder continuous. Moreover, they are real analytic in the open set where their gradients are not zero. No such criterion for real analyticity is valid  for the $\infty$-harmonic functions. One reason is that if $\nabla u = 0$ at some point,
then one may add an extra variable $x_{n+1}$ so that
$\nabla v \not = 0$ for 
$$v(x_1,x_2,\dots,x_n,x_{n+1}) = u(x_1,x_2,\dots,x_n) + x_{n+1}.$$
Both  are $\infty$-harmonic functions. Thus one always reaches the local situation with merely non-critical points.

\section{Viscosity Solutions}
\label{secvis}

\texttt{Jag undrar, sa flundran, om g\"{a}ddan \"{a}r fisk.}

\bigskip

The explicit presence of second derivatives in the $\infty$-Laplace Equation leads to difficulties, because they do not always exist when they are called for. This lack is the source of many problems. No weak formulation involving only first derivatives seems to be possible. It was first observed in [BdM] that the use of viscosity solutions circumvents this problem. First, let us explain why one does not restrict oneself only to $C^2$-solutions. The reason is that in the Dirichlet problem one can prescribe smooth boundary values so that no  $C^2$-solution can attain them. The critical points ($\nabla u = 0$) do not allow smoothness!

\begin{thm}[Aronsson] Suppose that $u \in C^2(\Omega)$ and that $\Delta_{\infty}u = 0$ in $\Omega.$ Then, either $\nabla u \not = 0$ in $\Omega$ or $u$ reduces to a constant.
\end{thm}

\emph{Proof:} The case $n = 2$ was proved by Aronsson in [A2], who used methods in Complex Analysis. The case $n\geq3$ is in [Y]. We skip the proof, since the theorem only serves us as a motivation.

\bigskip

The theorem has a  fatal consequence for the Dirichlet problem
\begin{equation}
\begin{cases} \Delta_{\infty}u = 0\quad\text{in}\quad \Omega\\
\phantom{ \Delta_{\infty}}u    = g\quad\text{on}\quad \partial \Omega
\end{cases}
\end{equation}
if $g \in C(\overline{\Omega})\cap C^2(\Omega)$ is deliberately chosen so that \emph{every} function $f 
\in C(\overline{\Omega})\cap C^2(\Omega)$ with boundary values $f = g$ on $\partial \Omega$ must have at least one critical point in $\Omega.$

 {\small For example, if $\Omega$ is the unit disc $x^2+y^2 < 1$ and $g(x,y) = xy,$ then every  $f$ with boundary values
$$f(\xi,\eta) = \xi\eta \quad \text{when} \quad \xi^2 + \eta^2 = 1,$$
be it a solution or not, must have a critical point in the open unit disc. This is a topological phenomenon related to index theory. Especially, a solution to the Dirichlet problem cannot be of class $C^2$ if such boundary values are prescribed, because it would contradict Aronsson's theorem. ---Another  way to see this is to use the symmetry. To this end, assume that $\nabla u \not = 0$ in the disc. By the end of Section \ref{secvar} the solution is unique. We can conclude that $u(x,y) = u(-x,-y).$ But upon differentiation we see that $\nabla u(0,0) = 0.$ There is a critical point in any case!}

The remedy is to introduce the concept of viscosity solutions. The modern theory of  viscosity solutions was developed by Crandall, Evans, Lions, Ishii, Jensen,and others. First, it was designed only for first order equations. Later, it was extended to second order equations. For the $\infty$-Laplacian it appeared in [BdM].

 {\small The label \emph{viscosity} came from first order equations like Burgers's Equation
\begin{equation}
\label{Burgers}
\begin{cases}
$$\dfrac{\partial u}{\partial t} + \dfrac{1}{2}\dfrac{\partial u^2}{\partial x} = 0\qquad(-\infty<x<\infty,\,t>0)$$\\
$$u(x,0) = g(x)$$
\end{cases}
\end{equation}
to which one adds a very small term representing viscosity. It is well known that even smooth initial values $g(x)$ to this Cauchy problem may force the solution $u = u(x,t)$ to develop shocks. Moreover, uniqueness often fails. A device was to add a term like $\varepsilon \Delta u,$ representing artificial viscosity: in this case
$$
\begin{cases}
\dfrac{\partial u_{\varepsilon}}{\partial t} + \dfrac{1}{2}\dfrac{\partial u^2_{\varepsilon}}{\partial x} =\varepsilon \dfrac{\partial^2 u_{\varepsilon}}{\partial x^2}\\
u_{\varepsilon}(x,0) = g(x)
\end{cases}$$
Now the solution $u_{\varepsilon} = u_{\varepsilon}(x,t)$ is smooth and unique. Moreover, $u_{\varepsilon}(x,t) \longrightarrow u(x,t) =$ some function. This $u$ is the desired viscosity solution of the original equation (\ref{Burgers}). This is the method of vanishing viscosity. However, just to identify the viscosity solution, the limit procedure with $\varepsilon$ can be skipped and the concept can be formulated directly for the original equation. This approach will be followed here.}

Let us return to the $\infty$-Laplace equation. Instead of considering the limit procedure
$$\Delta _{\infty}u_{\varepsilon} + \varepsilon \Delta u_{\varepsilon}\,=\,0$$
as $\varepsilon \to 0,$ we formulate the concept of viscosity solutions directly for the original equation\footnote{In the proof of everywhere differentiability in [ES] the approach with the vanishing term $\varepsilon \Delta u_{\varepsilon}$ seems to be essential.}. The above equation happens to be the Euler-Lagrange equation for the variational integral
$$J(v) = \int_{\Omega}\!e^{\frac{|\nabla v|^2}{2\varepsilon}}dx.$$

We begin with some infinitesimal calculus. Suppose that $v:\Omega \longrightarrow \R$ is a given function. Assume that $\phi \in C^2(\Omega)$ is \emph{touching $v$ from below} at some point $x_0 \in \Omega:$
\begin{equation}
\label{strict}
\begin{cases}
v(x_0) = \phi(x_0)\\
v(x) > \phi(x),\quad \text{when}\quad x \not = x_0
\end{cases}
\end{equation}
If it so happens that $v$ is smooth, then
$$
\begin{cases}
\nabla v(x_0) = \nabla  \phi(x_0)\\
\mathrm{D}^2v(x_0) \geq \mathrm{D}^2\phi(x_0)
\end{cases}$$
 by the infinitesimal calculus. Here
$$\mathrm{D}^2v(x_0) = \left(\frac{\partial^2v}{\partial x_i \partial x_j}(x_0)\!\right)_{\!\!n\times n}$$
is the Hessian matrix evaluated at the touching point $x_0.$ (For a symmetric matrix $\mathrm{A} = (a_{ij})$ we use the ordering
$$ \mathrm{A} \geq 0\quad \Longleftrightarrow\quad \sum_{1,j=1}^{n}a_{ij}\xi_i\xi_j \geq 0 \quad \text{for all}\quad \xi = (\xi_1,\xi_2,\dots,\xi_n).$$
Above $ \mathrm{A} = \mathrm{D}^2(v-\phi)(x_0).$) In particular, it follows that $v_{x_jx_j}(x_0) \geq 
\phi_{x_jx_j}(x_0),\, j = 1,2,\dots,n,$ and hence
\begin{align*}
\Delta v(x_0)& \,\geq\, \Delta \phi(x_0),\\
\Delta_{\infty} v(x_0)& \,\geq\, \Delta_{\infty} \phi(x_0),\\
\Delta_{p} v(x_0)& \,\geq\, \Delta_{p} \phi(x_0).\\
\end{align*}
Keep $n < p \leq \infty.$ If $v$ is a supersolution, i.e., $\Delta_{p} v \leq 0,$ then
$$ \Delta_{p} \phi(x_0) \leq 0.$$
The last inequality makes sense also if $v$ does not have any derivatives. This observation is important for the next definition.

\begin{definition}
\label{viscositydef}
 Let $n < p \leq \infty.$ We say that $v \in C(\Omega)$ is a \emph{viscosity supersolution} of the equation $\Delta_{p} v = 0$ in $\Omega,$ if
$$ \Delta_{p} \phi(x_0) \leq 0$$
 whenever $x_0 \in \Omega$ and  $\phi \in C^2(\Omega)$  are such that $\phi$ touches $v$ from below at $x_0.$ 

We say that $u \in C(\Omega)$ is a \emph{viscosity subsolution} of the equation $\Delta_{p} u = 0$ in $\Omega,$ if
$$ \Delta_{p} \psi(x_0) \geq 0$$
 whenever $x_0 \in \Omega$ and  $\psi \in C^2(\Omega)$  are such that $\psi$ touches $u$ from above at $x_0.$ 

Finally, we say that $h \in C(\Omega)$ is a \emph{viscosity solution}, if it is both a viscosity supersolution and a viscosity subsolution.
\end{definition}

\bigskip

Several comments are appropriate.
\begin{itemize}
\item It is often convenient to replace the condition of touching from below by requiring that
$$v-\phi\quad \text{has a (strict) minimum at}\quad x_0.$$
This gives an \textsf{equivalent definition}.
\item The differential operator is evaluated only at the touching point $x_0.$
\item Each point has its own family of test functions (it may be empty).
\item If there does not exist any function in $C^2(\Omega)$ touching at the point $x_0,$
 then there is no requirement. The point passes for free.
\item We required \emph{strict} inequality when $x\not = x_0.$ But using for example
$$\phi(x) + |x-x_0|^4$$
one can reach this from the requirement $\phi(x) \leq v(x).$ In several proofs \textsf{it is later needed that the touching is strict}.
\item Actually, it is possible to restrict the class of test functions to the second order polynomials
$$\phi(x) = c + \langle b,x \rangle + \frac{1}{2}\langle \mathrm{A}x,x \rangle.$$
\end{itemize}

\medskip

\begin{prop} [Consistency] A function $u \in C^2(\Omega)$ is a viscosity solution of $\Delta_pu = 0$ if and only if  $\Delta_pu(x) = 0$ holds pointwise in $\Omega.$
\end{prop}

\bigskip

\emph{Proof:} Let us consider the case of subsolutions. 

If $u$ is a viscosity subsolution, then $u$ itself will do as testfunction at the point $x_0,$ so that $\Delta_pu(x_0) \geq 0.$ (To write this in accordance with formula  (\ref{strict}) one can use the test function
$$ \psi(x) = u(x) + |x-x_0|^4,$$
so that the touching from above is \emph{strict} when $x \not = x_0.$ Then $0 \leq \Delta_p\psi(x_0) = \Delta_pu(x_0) .$)

If $u$ is a classical subsolution, i.\,e. $\Delta_pu(x)\geq 0$ at every point, we have to verify Definition \ref{viscositydef}. For a  test function  touching from above at $x_0$  the inequality
$$ \Delta_p\psi(x_0) \geq \Delta_pu(x_0)$$
holds
by the infinitesimal calculus. Since  $\Delta_pu(x_0)\geq 0$, the required inequality $ \Delta_p\psi(x_0) \geq 0$ follows. $\Box$

\bigskip

\emph{Example:} The function
$$u(x,y) = x^{\frac{4}{3}} - y^{\frac{4}{3}}$$
is a viscosity solution of $\Delta_{\infty}u = 0$ in the whole $\R^2$. To verify this is a direct calculation. The second derivatives blow up at the coordinate axes. For example, at the origin there does not exist any proper test function. We leave the verification to the reader. ---In fact,
$$u \in C_{\
loc}^{1,\frac{1}{3}}(\R^2) \cap W^{2,\frac{2}{3}-\varepsilon}_{loc}(\R^2),\qquad \varepsilon > 0.$$ {\small It has been conjectured that every viscosity solution in the plane belongs to $C_{\
loc}^{1,\frac{1}{3}}(\R^2)$, i.e. its gradient is locally H\"{o}lder continuous with exponent $ \tfrac{1}{3}.$ See [S] and [ES] for the  H\"{o}lder continuity.}

\bigskip

\emph{Example:} The function $|x-x_0|$ is a viscosity subsolution in the whole $\Rn.$ {\small Indeed, even the
function defined by the superposition\footnote{Needless to say, the Principle of Superposition is not generally valid for our nonlinear equation!}
$$V(x) = \int_{\Rn}\!|x-y|\rho(y)\,dy\qquad\text{where}\qquad \rho \geq 0$$
is a viscosity subsolution, if $ \int_{\Rn}\!|y|\rho(y)\,dy < \infty.$ It is a curious coincidence that the fundamental solution of the \emph{biharmonic} equation $\Delta\Delta u = - \delta$ in three dimensional space is $\tfrac{|x-x_0|}{8\pi}$ so that $\Delta\Delta V(x) = - 8\pi\rho(x).$ Given $V$, this tells us how to find a suitable $\rho$!}

\bigskip

In Section \ref{secvar} we have constructed the \emph{variational} solution
$$u_{\infty} = \lim_{j \to \infty} u_{p_j}$$
via a sequence of solutions to $p$-Laplace equations. In order to prove that the obtained $u_{\infty}$ is a \emph{viscosity} solution we must first establish that $u_{p_j}$ is a viscosity solution of the $p_{j}$-Laplace equation.

\begin{lemma}
\label{vill}
Suppose that $v \in C(\Omega)\cap W^{1,p}(\Omega)$ satisfies $\Delta_{p}v\leq 0$ in the weak sense:
$$\int_{\Omega}\!\langle|\nabla v|^{p-2}\nabla v,\nabla \eta\rangle\,dx\, \geq\, 0$$
for all $\eta \geq 0,\, \eta \in C_0^{\infty}(\Omega).$ Then $v$ is also a viscosity supersolution in $\Omega.$
\end{lemma}

\bigskip

\emph{Proof:} \textsl{Reductio ad absurdum.} If the claim is not valid, we have at some point $x_0$ in $\Omega$  a test function $\phi$ touching from below  so that
the inequality
$$\Delta_p \phi(x_0) > 0$$
holds. By continuity
$$\Delta_p \phi(x) > 0 \quad \text{when} \quad |x-x_0| < 2\rho$$
for some small radius $\rho > 0.$ Thus $\phi$ is a classical subsolution in $B(x_0,2\rho).$ Define
$$\psi(x) = \phi(x) + \frac{1}{2}\min_{\partial B(x_0,\rho)}\{v-\phi\} =  \phi(x) + \frac{1}{2}m.$$
Clearly, $\psi < v$ on $\partial B(x_0,\rho)$ and $\psi(x_0) > v(x_0).$ Let $D_{\rho}$ denote that component of the open set $\{\psi>v\}\cap B(x_0,\rho)$ which contains the point $x_0.$ We can conclude that $\psi = v$ on $\partial D_{\rho}$. 

Since $\psi$ is a weak subsolution and $v$ a weak supersolution of the same equation, we must by the Comparison Principle\footnote{The Comparison Principle is easy to verify. Insert the test function $\eta = [\psi-v]_+$ into 
$$\int_{D_{\rho}}\!\langle| \nabla \psi|^{p-2}\nabla \psi,\nabla \eta\rangle\,dx \leq 0, \quad 
 \int_{D_{\rho}}\!\langle| \nabla v|^{p-2}\nabla v,\nabla \eta\rangle\,dx \geq 0$$
and subtract to arrive at
$$\int_{D_{\rho}}\!\langle| \nabla \psi|^{p-2}\nabla \psi - | \nabla v|^{p-2}\nabla v,\nabla (\psi-v)_{+}\rangle\,dx \leq 0$$
By an elementary inequality, the integral is greater than or equal to  
$$2^{p-2}\int_{D_{\rho}}\!|\nabla(\psi -v)_+|^p\,dx.$$ Hence $\nabla(\psi -v)_+ =0.$ The result follows.}
 have $\psi \leq v$ in $D_{\rho}.$ This leads to the contradiction $\psi(x_0) \leq v(x_0) < \psi(x_0).$ 
The lemma follows.\qquad $\Box$

\bigskip

Then we isolate a fact from calculus.

\begin{lemma}
\label{convergence}
 Suppose that $f_j \longrightarrow f$ uniformly in $\overline{\Omega},$\, $f_j,f \in C(\overline{\Omega)}.$ If $\phi$ touches $f$ from below at $x_0 \in \Omega,$ then there are points $x_j \longrightarrow x_0$ such that
$$f_{j}(x_j) -\phi(x_j)\,=\, \min_{\Omega}\{f_{j}-\phi\}$$
for some subsequence.
\end{lemma}

\emph{Proof:} Since
$f_j-\phi = (f-\phi) + (f_j-f)$ and $f_j \longrightarrow f,$
$$\inf_{\Omega \setminus B_r}\{f_j-\phi\}\, \geq \,\tfrac{1}{2}\inf_{\Omega \setminus B_r}\{f-\phi\}\, >\, 0$$
when $j$ is large enough. Here $B_r = B(x_0,r)$ is a small ball. (It was essential to have the \emph{strict} inequality $\phi(x) < f(x)$ when $x \not = x_0.$) Hence
$$\inf_{\Omega \setminus B_r}\{f_j-\phi\} > f_j(x_0)- \phi(x_0)$$
when $j > j_r.$ This is so, because the right-hand side approaches zero while the left-hand side has a positive limit. Thus the value at the center is smaller than the infimum over $\Omega \setminus B_r.$ Therefore there is a point $x_j \in \overline{B(x_0,r)}$ such that
$$\min_{\Omega}\{f_j-\phi\} = f_j(x_j)-\phi(x_j)$$
when $j > j_r.$ To finish the proof, let $r \to 0$ via a sequence, say $r = 1,\tfrac{1}{2},\tfrac{1}{3},\dots.$\qquad   $\Box$

\bigskip

We aim at showing that the variational solution constructed in Section \ref{secvar} is a viscosity solution. We consider the slightly more general equation
$$\Delta_{p}v = - \varepsilon^{p-1},$$
which will be needed later. We repeat the variational procedure.

Again, let $f \in C(\overline{\Omega})$ be a Lipschitz continuous function with $\|\nabla f\|_{L^{\infty}(\Omega)} \leq L.$ Fix a number $\varepsilon \geqq 0$ and consider the problem of minimizing the variational integral
 $$J_p(v) = \int_{\Omega}\!\Bigl(\frac{1}{p}|\nabla v|^p - \varepsilon^{p-1}v\Bigr) dx,\qquad p > n,$$
among all functions in the class $ C(\overline{\Omega}) \cap W^{1,p}(\Omega)$ having boundary values $v = f$ on $\partial \Omega.$ The direct method in the Calculus of Variations shows the existence of a unique minimizer $v_p \in C(\overline{\Omega}) \cap W^{1,p}(\Omega),\,v_p|_{\partial \Omega} = f.$ 

We need a uniform bound in $p$ of the norms $\|\nabla v_p\|_{L^p(\Omega)}.$ Since $J_p(v_p) \leq J_p(f),$ we have
$$\frac{1}{p}\int_{\Omega}\!|\nabla v_p|^p\,dx \leq \frac{1}{p}\int_{\Omega}\!|\nabla f|^p\,dx + \varepsilon^{p-1}\int_{\Omega}\!(v_p -f)\,dx$$
and using Young's inequality
\begin{align*}
&\Bigl|\int_{\Omega}\!(v_p -f)\,dx\Bigr|\\
 \leq&\, \frac{\lambda^{p}}{p}\int_{\Omega}\!|v_p -f|^p\,dx + \frac{\lambda^{-q}}{q}|\Omega|\\
\leq&\, \frac{\lambda^{p}}{p}(\mathrm{diam}\Omega)^p\int_{\Omega}\!|\nabla v_p - \nabla f|^p \,dx +  \frac{\lambda^{-q}}{q}|\Omega|\\
\leq&\, \frac{2\lambda^{p}}{p}(\mathrm{diam}\Omega)^p\biggl\{\int_{\Omega}\!|\nabla v_p|^p\,dx +\int_{\Omega}\!|\nabla f|^p \,dx\biggr\} +  \frac{\lambda^{-q}}{q}|\Omega|,
\end{align*}
where a version of Friedrichs's inequality was used, and $q =\frac{p}{p-1}$ . (Needless to say, a more elegant estimation is available.) Selecting the auxiliary parameter $\lambda$ so that
for example
$\varepsilon^{p-1}\lambda^p(2\,\mathrm{diam}\Omega)^p = \tfrac{1}{2},$ we can absorb a term:
$$\frac{1}{2p}\int_{\Omega}\! |\nabla v_p|^p\,dx  \leq  \frac{3}{2p}\int_{\Omega}\! |\nabla f|^p\,dx 
+ \frac{\lambda^{-q}}{q}|\Omega|\varepsilon^{p-1}.$$
It follows that
$$\limsup_{p \to \infty}\|\nabla v_p\|_{L^p(\Omega)} \,\leq \, L + \varepsilon.$$
We also know that
\begin{align*}
|v_p(x) -v_p(y)| \leq\, &|\bigl(v_p(x)\! -\! f(x)\bigr) - \bigl(v_p(y)\!-\!f(y)\bigr)| + |f(x)\! -\!f(y)|\\
\leq &\,C_{\Omega}|x-y|^{1-\frac{n}{p}}\|\nabla (v_p-f)\|_{L^p(\Omega)} + L|x-y|
\end{align*}
from which we can see that the family $\{v_p\}$ is equicontinuous. By Ascoli's theorem we can extract a sequence $p_1,p_2,\dots$ such that
$$\begin{cases}
v_{p_j} \longrightarrow v_{\infty}\quad \text{uniformly in} \quad \overline{\Omega}\\
$$\nabla v_{p_j} \longrightarrow \nabla v_{\infty}\quad \text{weakly in every} \quad L^s(\Omega),\quad s<\infty,
\end{cases}$$
where $v_{\infty} \in C(\overline{\Omega}) \cap W^{1,\infty}(\Omega)$ is some function. That the same sequence will do simultaneously in each $L^s(\Omega)$ requires a diagonalization procedure. 

\bigskip

\begin{thm}
\label{existence}
 The function $ v_{\infty}$ is a viscosity solution of the equation
$$\max\{\varepsilon - |\nabla  v_{\infty}|,\, \Delta_{\infty} v_{\infty} \}\,=\,0.$$
In other words,
$$\varepsilon \leq |\nabla \phi(x_0)|\quad \text{and} \quad  \Delta_{\infty} \phi(x_0)\leq 0$$
holds for a test function touching from below and
$$\varepsilon \geq |\nabla \psi(x_0)|\quad \text{or}\quad  \Delta_{\infty} \psi(x_0)\geq 0$$
holds for test functions touching from above.

In the special case $\varepsilon = 0$ the equation reads $ \Delta_{\infty}v_{\infty} = 0.$
\end{thm}

\bigskip

\emph{Proof:} To prove the theorem, we first verify that $v_p$ is a viscosity solution of its own Euler-Lagrange equation
\begin{equation}
\label{Poisson}
\Delta_pv = - \varepsilon^{p-1}.
\end{equation}
In weak form, it reads 
$$\int_{\Omega}\!\langle |\nabla v|^{p-2}\nabla v,\nabla \eta\rangle\,dx\,=\,\varepsilon^{p-1}\int_{\Omega}\!\eta\,dx,$$
where $\eta \in C_0^{\infty}(\Omega).$ If $v$ were of class $C^2(\Omega)$ (and $p \geq 2$), this would be equivalent to the pointwise equation (\ref{Poisson}). The proof is the same as in Lemma \ref{vill}.

 Therefore, let us proceed with $p_j \to \infty.$ In the case of supersolutions we have that
$$\Delta_{p_j}\varphi(x_0) \leq -\varepsilon^{p-1}$$
with $\varphi$ touching $v_{p_j}$ from below at the point $x_0.$ Suppose that $\phi$ touches $v_{\infty}$ from below at $x_0,$ i.\,e. $\phi(x_0) = v_{\infty}(x_0)$ and $\phi(x) < v_{\infty}(x)$ when $x \not = x_0.$ By Lemma \ref{convergence} there are points $x_k \to x_0$ so that $v_{p_{j_k}} - \phi$ attains its minimum at $x_k.$ Thus 
$$\Delta_{p_{j_k}}\phi(x_k) \leq - \varepsilon^{p-1}$$
 or written out explicitely
$$|\nabla \phi(x_k)|^{p_{j_k}-4}\Bigl\{|\nabla \phi(x_k)|^{2}\Delta \phi(x_k) + (p_{j_k}-2)\Delta_{\infty}\phi(x_k)\Bigr\} \leq -\varepsilon^{p_{j_k}-1}.$$
If $\varepsilon \not = 0,$ it is impossible that $\nabla \phi(x_k) = 0,$ and we can divide out 
to obtain
\begin{equation}
\label{above}
 |\nabla \phi(x_k)|^{2}\frac{\Delta \phi(x_k)}{p_{j_k}-2} + \Delta_{\infty}\phi(x_k) \leq
-\frac{\varepsilon^3}{p_{j_k}-2}\Bigl(\frac{\varepsilon}{|\nabla \phi(x_k)|}\Bigr)^{p_{j_k}-4}.
\end{equation}
By continuity, the left-hand side has the limit $\Delta_{\infty} \phi(x_0).$ If $|\nabla \phi(x_0)| < \varepsilon,$ we get the contradiction $\Delta_{\infty} \phi(x_0) = - \infty.$ Therefore we  have established the first required inequality: $\varepsilon - |\nabla \phi(x_0)| \leq 0.$ This inequality implies that the last term in (\ref{above}) approaches $0$. Hence also
the second required inequality is valid: $\Delta_{\infty}\phi(x_0) \leq 0.$ This concludes the case $\varepsilon > 0.$

If $\varepsilon = 0$ there is nothing to prove, when $\nabla \phi(x_0) = 0.$ If $\nabla \phi(x_0) \not = 0$  then 
$$\frac{\Delta \phi(x_k)}{p_{j_k}-2} + \frac{\Delta_{\infty}  \phi(x_k)} {|\nabla \phi(x_k)|^2}  \leq 0$$
for large indices $k$. Again the desired result $\Delta_{\infty} \phi(x_0)
\leq 0$ follows. (In this case the first order equation is void.)
This concludes the case with supersolutions.

The case of subsolutions is even simpler and we omit it.\qquad $\Box$

\section{An Asymptotic Mean Value Formula}
\label{secasy}

The $\infty-$harmonic functions enjoy a kind of ''mean value'' property. The formula was suspected via considerations in game theory. To place it in a proper context, recall that a continuous function $u$ is harmonic if and only if it obeys the mean value formula discovered by Gauss. It is less well known\footnote{\textsc{W. Blaschke.} {\it Ein Mittelwertsatz und eine kennzeichende Eigenschaft des logarithmischen Potentials}, Leipz. Ber. \textbf{68}, 1916, pp. 3--7. and \textsc{I. Privaloff.} {\it Sur les fonctions harmoniques, Recueil Math\'{e}matique. Moscou (Mat. Sbornik)} \textbf{32}, 1925, pp. 464--471.} that $\Delta u = 0$ in $\Omega$ if and only if
\begin{equation}
 \label{as}
u(x) = \vint_{B(x,\varepsilon)}\!u(y)\,dy + o(\varepsilon^2) \qquad \text{as}\qquad \varepsilon \to 0 
\end{equation}
at every point $x \in \Omega.$ (We deliberately ignore the fact that the error term $o(\varepsilon^2) \equiv 0.$) We have taken the solid mean value. In [MPR] it was discovered that the $p$-harmonic functions have a counterpart, a remarkable mean value formula of their own:
\begin{equation*}
u(x) = \frac{p-2}{p+n}\,\frac{\underset{\overline{B(x,\varepsilon)}}{\max}\{u\} + \underset{\overline{B(x,\varepsilon)}}{\min}\{u\}}{2} +
\frac{2+n}{p+n}\vint_{B(x,\varepsilon)}\!u(y)\,dy +  o(\varepsilon^2) 
\end{equation*}
as $\varepsilon \to 0.$ The second term is a linear one, while the first term counts for the non-linearity: the greater the $p$, the more non-linear the formula. It has an interpretation in game theory, where it is essential that the coefficients (now probabilities!) add up to $1$:
$$ \frac{p-2}{p+n} + \frac{2+n}{p+n} = 1.$$ The cases $p =2$ and $p = \infty$ are noteworthy. When $p = 2$ the formula reduces to equation (\ref{as}). When $p=\infty$ we get the asymptotic formula

\medskip

\begin{equation}
\label{asym}
\boxed{u(x)\, =\, \dfrac{\underset{\overline{B(x,\varepsilon)}}{\max}\{u\} + \underset{\overline{B(x,\varepsilon)}}{\min}\{u\}}{2}  + 
o(\varepsilon^2) \qquad \text{as}\qquad \varepsilon \to 0 .}
\end{equation}

\medskip

This important formula provides an alternative definition. It turns out that a function $u\in C(\Omega)$ is a viscosity solution of $\Delta_{\infty}u = 0$ in $\Omega$ if and only if the above equation holds at each point in $\Omega,$ though interpreted in a specific way to be described below. Indeed, formula (\ref{asym}) is not valid directly for $u$ as it is written.
A counter example is
$$u(x,y) = x^{\frac{4}{3}} - y^{\frac{4}{3}}.$$
This $\infty$-harmonic function does not satisfy  (\ref{asym}) for example at the point $(1,0),$ which an exact computation shows, see [MPR]. 

The basic fact on which this depends is a property for test functions. The critical points in the domain are avoided.

\begin{lemma} Suppose that $\phi \in C^2(\Omega)$ and that $\nabla \phi(x_0) \not = 0$ at the point $x_0 \in \Omega.$ Then the asymptotic formula
\begin{equation}
\label{Taylor}
\boxed{\phi(x_0)\,=\,  \dfrac{\underset{\overline{B(x_0,\varepsilon)}}{\max}\{\phi\} + \underset{\overline{B(x_0,\varepsilon)}}{\min}\{\phi\}}{2}  + \dfrac{\Delta_{\infty}\phi(x_0)}{|\nabla \phi(x_0)|^2} +
o(\varepsilon^2) \quad \text{as}\quad \varepsilon \to 0}
\end{equation}
holds.
\end{lemma}

\bigskip

\emph{Proof:} We use Taylor's formula and some infinitesimal calculus to verify this. Since $\nabla \phi(x_0) \not = 0$ we can take $\varepsilon$ so small that $|\nabla \phi(x)| > 0$ when
$|x-x_0| \leq \varepsilon.$ The extremal values in {\small $\overline{B(x_0,\varepsilon)}$ }are attained  at points for which
\begin{equation}
\label{gra}
 x = x_0 \pm \varepsilon \frac{\nabla \phi(x)}{|\nabla \phi(x)|}, \qquad |x-x)| =\varepsilon.
\end{equation}
To see this, use for example a Lagrangian multiplier. Moreover, the $+$ -sign corresponds to the maximum, when $\varepsilon$ is small. The maximum and minimum points are \emph{approximately} opposite endpoints of some diameter. One also has
\begin{equation}
\label{epsgrad}
\frac{\nabla \phi(x)}{|\nabla \phi(x)|} = \frac{\nabla \phi(x_0)}{|\nabla \phi(x_0)|} + O(\varepsilon)\qquad \text{as} \qquad \varepsilon \to 0.
\end{equation}
For $x\in \partial B(x_0,\varepsilon)$ let $x^{\star}$ denote  \emph{exactly} opposite endpoints of a diameter, i.e., $x + x^{\star} = 2x_0.$ Adding the Taylor expansions
$$\phi(y) = \phi(x_0) + \big\langle \nabla \phi(x_0),y-x_0\big\rangle + \frac{1}{2}\big\langle \mathbb{D}^2\phi(x_0)(y-x_0),y-x_0\big\rangle + o(|y-x_0|^2)$$
for $y = x$ and $y = x^{\star}$, we can write
$$\phi(x) + \phi(x^{\star}) = 2 \phi(x_0) + \big\langle \mathbb{D}^2\phi(x_0)(x-x_0),x-x_0\big\rangle + o(|x-x_0|^2).$$
Note that the first order terms have cancelled.

Select $x$ to be the maximal point:
$$\phi(x) = \max_{|y-x_0|=\varepsilon}\{\phi(y)\}.$$
Inserting  formula (\ref{gra}) with the approximation (\ref{epsgrad}), we obtain the inequality
\begin{align*}
&\phantom{AB}\underset{\overline{B(x_0,\varepsilon)}}{\max}\{\phi\} + \underset{\overline{B(x_0,\varepsilon)}}{\min}\{\phi\} \leq\, \phi(x) + \phi(x^{\star})\\ &=\,
 2\, \phi(x_0) + \big\langle \mathbb{D}^2\phi(x_0)(x-x_0),x-x_0\big\rangle + o(|x-x_0|^2)\\ &=\,
 2\, \phi(x_0) + \Big\langle \mathbb{D}^2\phi(x_0)\frac{\nabla \phi(x_0)}{|\nabla \phi(x_0)|},\frac{\nabla \phi(x_0)}{|\nabla \phi(x_0)|}\Big\rangle + o(\varepsilon^2)\\&
= \, 2\, \phi(x_0) +\frac{\Delta_{\infty}\phi(x_0)}{|\nabla \phi(x_0)|^2} + o(\varepsilon^2).
\end{align*}

Then, selecting $x$ to be a minimal point:
$$\phi(x) = \min_{|y-x_0|=\varepsilon}\{\phi(y)\}$$
we can derive the opposite inequality, since now
$$\underset{\overline{B(x_0,\varepsilon)}}{\max}\{\phi\} + \underset{\overline{B(x_0,\varepsilon)}}{\min}\{\phi\}\, \geq \,\phi(x^{\star}) + \phi(x).$$
This proves the asymptotic formula (\ref{Taylor}).\qquad
$\Box$

\bigskip

There is a slight problem about formula (\ref{Taylor}) when $\nabla \phi(x_0) =  0.$
Therefore the definition below contains amendments for this situation, excluding some test functions.

\bigskip

\begin{definition}
\label{viscmean}
A function $u \in C(\Omega)$ satisfies the formula
\begin{equation}
\label{asymm}
u(x)\, =\, \dfrac{\underset{\overline{B(x,\varepsilon)}}{\max}\{u\} + \underset{\overline{B(x,\varepsilon)}}{\min}\{u\}}{2}  + 
o(\varepsilon^2) \qquad \text{as}\qquad \varepsilon \to 0 
\end{equation}
\emph{in the viscosity sense}, if the two conditions below hold:
\begin{itemize}
\item if  $x_0 \in \Omega$ and if $\phi \in C^2(\Omega)$ touches $u$ from below at $x_0,$ then
$$\phi(x_0) \geq \frac{1}{2}\Bigl(\underset{|y-x_0| \leq \varepsilon}{\max}\{\phi(y)\} + \underset{|y-x_0| \leq \varepsilon}{\min}\{\phi(y)\}\Bigr)  + 
o(\varepsilon^2) \qquad \text{as}\qquad \varepsilon \to 0 .$$
Moreover, if it so happens that $\nabla \phi(x_0) = 0,$ we require  that the test function satisfies $\mathbb{D}^2\phi(x_0) \leq 0.$
\item  if  $x_0 \in \Omega$ and if $\psi \in C^2(\Omega)$ touches $u$ from above at $x_0,$ then
$$\psi(x_0) \leq \frac{1}{2}\Bigl(\underset{|y-x_0| \leq \varepsilon}{\max}\{\psi(y)\} + \underset{|y-x_0| \leq \varepsilon}{\min}\{\psi(y)\}\Bigr)  + 
o(\varepsilon^2) \qquad \text{as}\qquad \varepsilon \to 0 .$$
Moreover, if it so happens that $\nabla \psi(x_0) = 0,$ we require  that the test function satisfies $\mathbb{D}^2\psi(x_0) \geq 0.$
\end{itemize}
\end{definition}

\bigskip

\emph{Remark:} The extra restrictions when $\nabla \phi(x_0) = 0$ or $\nabla \psi(x_0) = 0$
are used in the form
$$\lim_{y\to x_0}\frac{\phi(y)-\phi(x_0)}{|y-x_0|^2}\,\leq \,0,\qquad \lim_{y\to x_0}\frac{\psi(y)-\psi(x_0)}{|y-x_0|^2}\,\geq \,0.$$

\bigskip

The reader may wish to verify that the previously mentioned function $u(x,y) = x^{4/3}-y^{4/3}$ satisfies (\ref{asymm}) in the viscosity sense.

The main result of this section is the characterization below.

\begin{thm} Let $u \in C(\Omega).$ Then $\Delta_{\infty}u = 0$ in the viscosity sense if and only if the mean value formula (\ref{asymm}) holds in the viscosity sense.
\end{thm}

\bigskip

\emph{Proof:} Let us consider  subsolutions. 
Suppose that $\psi$ touches $u$ from above at the point $x_0 \in \Omega.$ If $\nabla \psi(x_0) \not = 0$ then, according to formula (\ref{Taylor}), the condition $\Delta_{\infty}\psi(x_0) \geq 0$ is equivalent 
to the inequality for $\psi$ in Definition \ref{viscmean}. Thus this case is clear.

Assume, therefore, that  $\nabla \psi(x_0)  = 0.$ First, if $u$ satisfies formula (\ref{asymm}) in the viscosity sense, the extra assumption
 $$\lim_{y\to x_0}\frac{\psi(y)-\psi(x_0)}{|y-x_0|^2}\,\geq \,0$$
is at our disposal.
Thus if
$$
\psi(x_{\varepsilon}) = \underset{|x-x_o|\leq \varepsilon}{\min}\psi(x)$$
we conclude that
\begin{align*}
&\,\liminf_{\varepsilon \to 0}\frac{1}{\varepsilon^2}\Bigl\{\frac{1}{2}\bigl(\max_{\overline{B(x_0,\varepsilon)}}\psi +
\min_{\overline{B(x_0,\varepsilon)}}\psi \bigr) - \psi(x_0)\Bigr\}\\
=
&\,\liminf_{\varepsilon \to 0}\frac{1}{\varepsilon^2}\Bigl\{\frac{1}{2}\bigl(\max_{\overline{B(x_0,\varepsilon)}}\psi -\psi(x_0)\bigr) +\frac{1}{2}\bigl(\min_{\overline{B(x_0,\varepsilon)}}\psi-\psi(x_0)\bigr) \Bigr\}\\
\geq&\,\liminf_{\varepsilon \to 0}\frac{1}{\varepsilon^2}\,\frac{1}{2}\bigl(\min_{\overline{B(x_0,\varepsilon)}}\psi -\psi(x_0)\bigr)\\
=&\,
\frac{1}{2}\liminf_{\varepsilon \to 0}\Bigl(\frac{\psi(x_{\varepsilon})-\psi(x_0)}{|x_{\varepsilon}-x_0|^2}\Bigr)\frac{|x_{\varepsilon}-x_0|^2}{\varepsilon^2}\,\geq\, 0
\end{align*}
since $\tfrac{|x_{\varepsilon}-x_0|^2}{\varepsilon^2} \leq 1.$ This proves the desired inequality in Definition \ref{viscmean} in the case $\nabla \phi(x_0)  = 0.$

On the other hand, when $\nabla \psi(x_0) = 0$   the condition $\Delta_{\infty}\psi(x_0) \geq 0$ is certainly fulfilled. Thus the mean value property in the viscosity sense guarantees that we have a viscosity subsolution of the $\infty$-Laplace equation.\qquad
$\Box$

\bigskip

We may add that the formula (\ref{asym}) is basic in image processing and for some numerical algorithms.

\section{Comparison with Cones}
\label{seccom}

The graph of the function 
$$C(x) = a + b|x-x_0|$$
is a half cone and $C =C(x)$ is a viscosity solution of the $\infty$-Laplace Equation, when $x\not=x_0.$ In the whole $\Rn$ it is a viscosity supersolution, if $b < 0.$  It plays the role of a fundamental solution of the $\infty$-Laplace equation and $|x-x_0|$ is the limit of the functions
$$|x-x_0|^{\frac{p-n}{p-1}}$$
encountered for the $p$-Laplace equation. The astonishing fact is that if a function obeys the comparison principle with respect to the cone functions, then it is $\infty$-harmonic and, consequently, it obeys the comparison principle    with respect to \emph{all} viscosity super- and subsolutions (proof later in Section \ref{secuni}). 

\begin{definition}
\label{conecomparison}
We say that that $u \in C(\Omega)$ obeys the comparison principle from above with cones if, for each subdomain $D \subset \subset \Omega$, the inequality
$$u(x)-b|x-x_0|\,\leq\, \max_{\xi \in \partial D}\bigl\{u(\xi)-b|\xi-x_0|\bigr\}$$
holds when $x \in D,\, x_0 \in \Rn \setminus D.$

 For the  comparison principle from below with cones the required inequality is
$$u(x)-b|x-x_0|\,\geq\, \min_{\xi \in \partial D}\bigl\{u(\xi)-b|\xi-x_0|\bigr\}.$$

The function obeys the comparison principle with cones, if it does so both from above and below.
\end{definition}

\bigskip

If the subdomain is a ball $B(x_0,r)\subset\subset\Omega$  we have the convenient inequality

\begin{equation}
\label{cone}
\boxed{
u(x)\, \leq\, u(x_0) + \underset{|\xi-x_0| = r}{\max}\Bigl(\dfrac{u(\xi)-u(x_0)}{r}\Bigr)|x-x_0|}
\end{equation}

\medskip 
\noindent when $|x-x_0| \leq r.$ This is so, because the inequality is always true at $x=x_0$ and when $x\in \partial B(x_0,r).$ (Take $D= \{x|\,0 < |x-x_0| < r\}.$) For cones the necessicity of the Comparison Principle is easy to prove\footnote{The general case is in Section \ref{secuni}.}:

\begin{prop}
\label{onsd}
 If $\Delta_{\infty}u \geq 0$ in the viscosity sense, the comparison principle with respect to cones from above holds.
\end{prop}

\bigskip

\emph{Proof:} We shall show that
\begin{equation}
\label{111}
u(x)-b|x-x_0|\,\leq\, \max_{\xi \in \partial D}\bigl\{u(\xi)-b|\xi-x_0|\bigr\},
\end{equation}
when $x \in D,\,x_0\not \in D.$
We first treat the case $b \not = 0.$ Now\footnote{For radial functions $\Delta_{\infty}R(r) = R'(r)^2R''(r),\,\, r =|x-x_0|.$}
\begin{align}
\label{222}
\Delta_{\infty}\bigl(b|x-x_0|-\gamma|x-x_0|^2\bigr)\\
= -2 \gamma\bigl(b-2\gamma|x-x_0|\bigr)^2\,<\,0,\nonumber
\end{align}
when the auxiliary parameter $\gamma > 0$ is small enough.

\textsf{Claim:}
$$u(x)-\overbrace{(b|x-x_0|-\gamma|x-x_0|^2)}^{\psi(x)} \leq \max_{\xi \in \partial D}\bigl\{u(\xi)  -\overbrace{(b|\xi-x_0|-\gamma|\xi-x_0|^2)}^{\psi(\xi)} \bigr\}.$$
The expression $u(x)-\psi(x)$ has a maximum in $\overline{D}.$ If it is attained in $D,$ then $\Delta_{\infty}\psi(x) \geq 0$ at that point, but this contradicts the above calculation in (\ref{222}). Hence the maximum is attained (only) on the boundary. Thus the claim holds. Now we let $\gamma \to 0$. The desired  inequality (\ref{111}) follows for the case $b \not = 0.$

The case $b=0$ follows immediately, by letting $b \to 0$ above in (\ref{111}). The inequality reduces to the maximum principle
$$u(x) \leq \max_{\xi \in \partial D}u(\xi).$$
As a byproduct, we got another proof of the \textsf{maximum principle.}\quad $\Box$

\bigskip

Comparison with cones yields a simple proof of Harnack's inequality, encountered in Section \ref{secvar}  for variational solutions.

\begin{prop} [Harnack]
\label{Har}
 Suppose that $\Delta_{\infty} v \leq 0$ in $\Omega$ in the viscosity sense. If $v \geq 0$ in the ball $B(x_0,R) \subset  \Omega,$ then
$$ v(y) \leq 3\,v(x)\quad \text{when}\quad x,y \in B(x_o,r),\,\,4r<R.$$
\end{prop}

\bigskip

\emph{Proof:} Comparison with cones in $D = B(y,3r) \setminus\{y\}$ yields
\begin{align*}
v(x) &\geq v(y) + \frac{|x-y|}{3r}\,\underset{|\xi-y| =3r}{\min}\{v(\xi)-v(y)\}\\
&= \Bigl(1-\frac{|x-y|}{3r}\Bigr)v(y) + \frac{|x-y|}{3r}\, \underset{|\xi-y| =3r}{\min}\{v(\xi)\}\\
& \geq \Bigl(1-\frac{|x-y|}{3r}\Bigr)v(y).
\end{align*}
The result follows, since $|x-y| < 2r.$ $\Box$

\bigskip

Comparison with cones immediately yields a counterpart to Hadamard's \textsf{Three Spheres Theorem}. If $\Delta_{\infty}u \geq 0$ in the annulus $r < |x| < R$ and if $u$ is continuous in the closure of the annulus, then
$$\boxed{u(x) \,\leq\, \dfrac{M_R\,(|x|-r)}{R-r} + \dfrac{M_r\,(R-|x|)}{R-r}\quad \text{where}\quad M_{\rho} = \max_{\partial B(0,\rho)}u,}$$
when $r \leq |x| \leq R.$ For viscosity supersolutions there is a similar inequality from below with minima on spheres. ---In other words:

\begin{prop} If $u$ is a viscosity subsolution to the $\infty-$ Laplace equation in the annulus
$r_1 < |x-x_0| < r_2,$ then the function
$$r \longmapsto \max_{\partial B(x_0,r)}u\qquad\text{is convex
when}\qquad r_1 < r < r_2.$$
For a viscosity supersolution $v$ the function
$$r \longmapsto \min_{\partial B(x_0,r)}v\qquad\text{is concave
when}\qquad r_1 < r < r_2.$$
\end{prop}

\bigskip

We can conclude that if an $\infty$-harmonic function is bounded around an isolated boundary point, then it has a limit at the point. (Yet, the singularity is not removable. Example: $|x|.$)

\begin{lemma}
\label{conetovisc}
If $u \in C(\Omega)$ obeys the comparison with cones from above, then $\Delta_{\infty}u \geq 0$ in the viscosity sense.
\end{lemma}

\bigskip

\emph{Proof:} Choose a point $x \in \Omega.$ Assume that the test function $\phi$ touches $u$ fram above at $x:$
$$\phi(x)-u(x) < \phi(y)-u(y) \quad \text{when} \quad y\not=x,\, y \,\in \Omega.$$
We have to show that $\Delta_{\infty}\phi(x) \geq 0.$ Denote $\mathfrak{p} = \nabla \phi(x).$ We first consider the case $\mathfrak{p} \not = 0. $ Let
$$y = x - \lambda \nabla \phi(x) =  x - \lambda \mathfrak{p}$$ be the center of the ball $B(y,r),$ where $\lambda$ is small (it will be sent to $0$). By the assumption
$$ \phi(x) \leq \phi(y) + |x-y|\underset{|\xi-y| =r}{\max}\frac{\phi(\xi)-\phi(y)}{r}$$
An elementary calculation transforms the inequality into 
\begin{equation}
\label{now}
\phi(x)-\phi(y) \leq \frac{|x-y|}{r-|x-y|}\,\underset{|\xi-y| =r}{\max}\left(\phi(\xi)-\phi(x)\right).
\end{equation}
Write $\mathbb{X} = D^2\phi(x)$ in Taylor's formula
\begin{align*}
\phi(y) = \phi(x) + \langle\mathfrak{p},y-x\rangle + \frac{1}{2}\langle \mathbb{X}(y-x),y-x\rangle + o(|y-x|^2)\\
= \phi(x) -\lambda |\mathfrak{p}|^2  + \frac{1}{2}\lambda^2\langle \mathbb{X}\mathfrak{p},\mathfrak{p}\rangle + \lambda^2 o(|\mathfrak{p}|^2)
\end{align*}
and use inequality (\ref{now}) to obtain
\begin{gather*}
|\mathfrak{p}|^2 -  \frac{1}{2}\lambda\langle \mathbb{X}\mathfrak{p},\mathfrak{p}\rangle 
+ \lambda o(|\mathfrak{p}|^2)\\
\leq \frac{|\mathfrak{p}|}{r-\lambda |\mathfrak{p}|}\, 
\underset{|\xi-y| =r}{\max}\left(\phi(\xi)-\phi(x)\right).
\end{gather*}
The maximum on the sphere $B(x,r)$ is attained at, say $\xi_r.$ Sending $\lambda$ to $0$ (and $y$ to $x$)
we obtain
$$|\mathfrak{p}|^2  \leq \frac{|\mathfrak{p}|}{r}\left(\phi(\xi_r)-\phi(x)\right),$$
where we can write, by Taylor's formula:
\begin{equation*}
\phi(\xi_r)-\phi(x) =\langle \mathfrak{p},\xi_r-x\rangle + \frac{1}{2}\langle \mathbb{X}(\xi_r-x),\xi_r-x\rangle + o(|\xi_r-x|^2).
\end{equation*}
It follows from Schwarz's Inequality that
$$|\mathfrak{p}|^2  \leq |\mathfrak{p}|^2 +\frac{|\mathfrak{p}|r}{2}\langle\mathbb{X}\,\frac{\xi_r-x}{r},\frac{\xi_r-x}{r}\rangle +|\mathfrak{p}| o(r) $$
and so
$$\langle\mathbb{X}\,\frac{\xi_r-x}{r},\frac{\xi_r-x}{r}\rangle +\frac{ o(r)}{r} \,\geq\, 0.$$
One can deduce that, at least via a sequence of $r$'s, 
$$\lim_{r\to 0}\frac{\xi_r-x}{r} = \frac{\mathfrak{p}}{|\mathfrak{p}|}.$$
In fact,  the same inequalities also imply that
$$|\mathfrak{p}| \leq \langle\mathfrak{p},\frac{\xi_r-x}{r}\rangle +O(r)   \leq \left|\langle\mathfrak{p},\frac{\xi_r-x}{r}\rangle\right| + O(r) \leq |\mathfrak{p}| + O(r),$$
from which the above limit is evident.

Thus we arrive at the inequality
$$\langle\mathbb{X}\mathfrak{p}, \mathfrak{p}\rangle\,\geq\,0.$$
We have assumed that $\mathfrak{p} \not = 0,$ but the inequality certainly holds without this restriction. This is the desired inequality $\Delta_{\infty}\phi(x) \geq 0.$ \qquad $\Box$

\bigskip

\begin{thm} A continuous function is a viscosity solution of the $\infty$-Laplace Equation if and only if it obeys the comparison with respect to cones.
\end{thm}

\bigskip

As a matter of fact, most of the theory of $\infty$-harmonic functions can be based on the Principle of Comparison with Cones, as in [U]. A remarkable simple uniqueness proof is in [AS], where the machinery with viscosity solutions is avoided.

\paragraph{Semicontinuous $u$ and $v$} In connexion with Perron's Method it is convenient to replace the continuity requirement in the definition for the viscosity supersolutions and subsolutions with semicontinuity. Thus we say that a function
$u: \Omega \to [-\infty,\infty)$ is a \emph{viscosity subsolution} of $\Delta_{\infty}u = 0,$
if 
\begin{description}
\item{(i)} $ u \not \equiv -\infty$,
\item{(ii)} $u$ is upper semicontinuous in $\Omega$,
\item{(iii)} whenever $\phi \in C^2(\Omega)$ touches $u$ from above at a point $x_0,$ we have $\Delta_{\infty}\phi(x_0) \geq 0.$
\end{description}
In other words, the continuity is weakened to conditions (i) and (ii), while the rest stays as before. ---The viscosity supersolutions are required to be lower semicontinuous and not identically $+\infty$. It is a routine to verify that \textsf{the results} of the present section \textsf{are valid under the new definition}. Actually, this seemingly weaker definition does not allow any ''new'' viscosity subsolutions. It follows from the next lemma that it is equivalent to the old definition.

\bigskip

\begin{lemma}
\label{semicont}
 The semicontinuous viscosity sub- and supersolutions are continuous.
\end{lemma}

\bigskip

\emph{Proof:} Consider the case of a semicontinuous viscosity subsolution $u$. It is easy to verify that Proposition \ref{onsd} holds with this change. So does Harnack's Inequality (Proposition \ref{Har}). Applying it to the function $v = -u$ we conclude, using a suitable chain of balls, that $u$ is \emph{locally bounded in} $\Omega.$ 

The upper semicontinuity is in the definition and so we have to prove only that $u$ is lower semicontinuous. To this end, let $x_0$ be an arbitrary point in $\Omega.$ Fix
$r>0$ so that $B(x_0,2r) \subset \subset \Omega$ and let $M = \max\{u\}$ taken over the closure of $B(x_0,2r).$ Select a sequence of points $x_j \in B(x_0,r)$ such that
$$x_j \to x_0\qquad \text{and} \qquad \lim_{j \to \infty} u(x_j) = \liminf_{x\to x_0}u(x).$$
Comparison with cones from above yields
\begin{align*}
u(x)\,& \leq\, u(x_j) + \max_{|\xi-x_j|=r}\Bigl(\frac{u(\xi)-u(x_0)}{r}\Bigr)|x-x_j|\\ \,&\leq
\,u(x_j) + \frac{2M}{r}|x-x_j|,
\end{align*}
when $|x-x_j| \leq r.$ Since $x_0 \in B(x_j,r)$ we see that
$$u(x_j) \geq u(x_0) - \frac{2M}{r}|x_0-x_j|$$
and hence $\lim u(x_j) \geq u(x_0).$ Thus $u$ is also lower semicontinuous.\qquad $\Box$

\bigskip

The same proof yields local Lipschitz continuity, see [J].

\section{From the Theory of Viscosity Solutions}
\label{secfro}

Convex and concave functions are important auxiliary tools. The main reason is that they have second derivatives of a kind based on the Taylor polynomial. Recall that a convex function  is locally Lipschitz continuous and, consequently, almost everywhere differentiable  according to Rademacher's theorem. The first derivatives are Sobolev derivatives. The second derivatives exist in the sense of Alexandrov.

\begin{thm} [Alexandrov] If  $v: \Rn \to \R$ is a convex function, then at almost every point $x$ there exists a symmetric $n\times n$-matrix $\mathbb{A} = \mathbb{A}(x)$ such that the expansion
$$v(y) = v(x) + \langle \nabla v(x),y-x\rangle + \frac{1}{2}\langle \mathbb{A}(x)(y-x),y-x\rangle + o(|x-y|^2)$$ 
is valid as   $y \to x.$
\end{thm}

For a proof we refer to [EG, Section 6.4, pp. 242--245]. When existing,
the second Sobolev derivatives $\mathbb{D}^2v$ of a convex function 
will do as Alexandrov derivatives. However, the Alexandrov derivatives are not always Sobolev derivatives: \emph{the formula for integration by parts may fail}. Nonetheless, the notation $\mathbb{A} = \mathbb{D}^2v$ is used. 

We seize the opportunity to define the
$$\text{\sf infimal convolution:}\qquad\quad 
v_{\varepsilon}(x) = \inf_{y\in\Omega}\Bigl\{v(y) + \frac{1}{2\varepsilon}|y-x|^2\Bigr\}$$
of a (lower semi)continuous function $v$, which for simplicity is assumed to be $\geq 0$ in $\Omega.$  At each point $x\in \Omega$
$$v_{\varepsilon}(x) \nearrow v(x).$$
An elementary calculation shows that the function
$$v_{\varepsilon}(x) - \dfrac{|x|^2}{2\varepsilon}$$
is  concave. Hence it has second Alexandrov derivatives and so has $v_{\varepsilon}.$ Thus the derivatives $\mathbb{D}^2v_{\varepsilon}$ are available (in the sense of Alexandrov). 
A remarkable property is that the infimal convolution preserves viscosity supersolutions of the $p$-Laplace equation, $p \leq \infty:$ 
$$\Delta_pv \geq 0\quad \Longrightarrow \quad \Delta_pv_{\varepsilon} \geq 0\qquad\text{(in the viscosity sense)}.$$

Equivalent definitions for viscosity solutions can be formulated in terms of so-called jets. Letters like $\mathbb{X}$ and $\mathbb{Y}$ denote symmetric real $n\times n$-matrices. For a function $v:\Omega \to \R$ we define the \emph{subjet} $J^{2,-}v(x)$ at the point $x \in \Omega$ as the set of all $(\mathfrak{p},\mathbb{X})$ satisfying the inequality 
$$v(y)\geq v(x) +\langle \mathfrak{p},y-x \rangle + \frac{1}{2}\langle \mathbb{X}(y-x),y-x\rangle + o(|x-y|^2)\quad\text{as} \quad y \to x,$$
where $\mathfrak{p}$ is a vector in $\Rn.$ Without the error term, the second order polynomial in the variable $y$ at the right-hand side would immediately do as a test function touching from below\footnote{As a mnemonic rule, \emph{sub}jets touch from below and are therefore used for testing \emph{super}solutions.} at the point $x$. Notice also that if $v\in C^2(\Omega)$, then one can use the Taylor polynomial with $\mathfrak{p} =\nabla v(x)$ and $\mathbb{X} = \mathrm{D}^2v(x).$

To define the \emph{superjet} $J^{2,+}u(x)$ we use the opposite inequality
$$u(y)\leq u(x) +\langle \mathfrak{p},y-x \rangle + \frac{1}{2}\langle \mathbb{X}(y-x),y-x\rangle + o(|x-y|^2)\quad\text{as} \quad y \to x.$$
---The Alexandrov derivatives are members of the jets. The jet at a given point may well be empty, but there are always points nearby at which the jet is non-empty; the simple proof is in [Ko, Proposition 2.5, p. 18].

\begin{prop} Let $v\in C(\Omega)$ and $x\in \Omega.$ There exist points $x_k \to x$ such that $J^{2,-}v(x_k)$ is non-empty when $k=1,2,\dots.$ The same holds for superjets.
\end{prop}

\bigskip

Finally, we introduce the \emph{closures of the jets}. The closure $\overline{J^{2,-}}v(x)$ contains all $(\mathfrak{p},\mathbb{X})$ such that there exists a sequence $$(x_k,v(x_k),\mathfrak{p}_k,\mathbb{X}_k) \longrightarrow (x,v(x),\mathfrak{p},\mathbb{X}),\,$$ where $x_k \in \Omega,\,\, (\mathfrak{p}_k,\mathbb{X}_k) \in J^{2,-}v(x_k).$
 
 In terms of jets the $p$-Laplace operator takes the form
\begin{equation}
F_p(\mathfrak{p},\mathbb{X}) =
\begin{cases}
|\mathfrak{p}|^{p-4}\left\{|\mathfrak{p}|^2\mathsf{Trace}(\mathbb{X}) + (p-2)\langle\mathbb{X}\mathfrak{p},\mathfrak{p}\rangle\right\},\quad p < \infty
\\\langle\mathbb{X}\mathfrak{p},\mathfrak{p}\rangle, \quad p = \infty
\end{cases}
\end{equation}

\bigskip

\begin{prop} For $v\in C(\Omega)$ the following conditions are equivalent:
\begin{itemize}
\item$ $v is a viscosity supersolution of $\Delta_pv \leq 0.$
\item $F_p(\mathfrak{p},\mathbb{X}) \leq 0$
when $(\mathfrak{p},\mathbb{X}) \in J^{2,-}v(x),\, x\in \Omega.$
\item $F_p(\mathfrak{p},\mathbb{X}) \leq 0$
when $(\mathfrak{p},\mathbb{X}) \in \overline{J^{2,-}}v(x),\, x\in \Omega.$
\end{itemize}
\end{prop}

\bigskip

It is enough that $v$ is lower semicontinuous. There is, of course, a similar characterization for viscosity subsolutions. This is Proposition 2.6 in [Ko], where a proof can be found. A special case of Ishii's Lemma (also called the Theorem of Sums) is given below. It requires doubling of the variables.  

\bigskip

\begin{thm}
\label{Ishii}
 [Ishii's Lemma] Let $u$ and $v$ belong to $C(\Omega).$ If there exists an interior point
$(x_j,y_j)$ in $\Omega \times \Omega$ for which the maximum 
$$\max_{x,y\in\Omega}\left\{u(x)-v(y) - \tfrac{j}{2}|x-y|^2\right\}$$
is attained,
 then there are symmetric matrices $\mathbb{X}_j$ and $\mathbb{Y}_j$ such that
$$\bigl(j(x_j-y_j),\mathbb{X}_j\bigr) \in \overline{J^{2,+}}u(x_j),\quad \bigl(j(x_j-y_j),\mathbb{Y}_j\bigr) \in \overline{J^{2,-}}v(y_j)$$
and $\mathbb{X}_j \leq \mathbb{Y}_j.$
\end{thm}

\bigskip

 Observe carefully that if it so happens that the involved functions are of class $C^2,$ then one would have
$$ j(x_j-y_j) = \nabla u(x_j) = \nabla v(y_j),\quad \mathrm{D}^2u(x_j) =  \mathbb{X}_j ,\quad    \mathrm{D}^2v(y_j) = \mathbb{Y}_j$$
according to the infinitesimal calculus. The lemma is valid also for semicontinuous functions. For the proof we refer to [Ko].

\section{Uniqueness of Viscosity Solutions}
\label{secuni}

\texttt{Unus, sed leo}

\bigskip
 
 We shall present the uniqueness proof, originally due to Jensen [J], for the Dirichlet problem in a bounded, but otherwise quite arbitrary domain. A different proof is in [BB]. A proof based on comparison with cones is in [AS].

\bigskip

\begin{thm}[Jensen]
\label{uniqueness}
Let $\Omega$ be an arbitrary bounded domain in $\Rn.$ Given a Lipschitz continuous function $f: \partial \Omega \to \R,$ there exists a unique viscosity solution $u\in C(\overline{\Omega})$ with boundary values $f.$ Moreover, $u \in W^{1,\infty}(\Omega),$ and $\|\nabla u\|_{L^{\infty}}$ has a bound depending only on the Lipschitz constant of $f.$
\end{thm}

\bigskip

The existence of a solution was proved in Theorem \ref{existence}. Thus we concentrate on the uniqueness. Following a device of Jensen we introduce two auxiliary equations with parameter $\varepsilon > 0.$ The situation will be
\begin{align*}
&\max\{\varepsilon -|\nabla u^+|,\,\Delta_{\infty}u^+\}\, =\,0\qquad\qquad&\text{Upper Equation}\\
&\Delta_{\infty}u\, =\, 0\qquad\qquad&\text{Equation}\\
&\min\{|\nabla u^-|-\varepsilon,\,\Delta_{\infty}u^-\}\,=\,0\qquad\qquad&\text{Lower Equation}
\end{align*}
As we shall see, when the functions have the same boundary values they are ordered: $u^- \leq u\leq u^+.$ We  need only viscosity supersolutions for the Upper Equation and only subsolutions for the Lower Equation. To recall the definition, see Theorem \ref{existence}. Given boundary values $f$, extended so that
$$f \in C(\overline{\Omega}),\quad \|\nabla f\|_{L^p(\Omega)} \leq L,$$ we constructed in Section \ref{secvis}  a viscosity solution to the Upper Equation via the minimizers of the variational integrals
$$J(v) = \int_{\Omega}\!\Bigl(\frac{1}{p}|\nabla v|^p - \varepsilon^{p-1}v \Bigl)\,dx.$$
We arrived at the bound below.
\bigskip

\begin{lemma}
\label{L}
 A variational solution of the Upper Equation satisfies
$$ \|\nabla v\|_{L^{\infty}(\Omega)} \leq L + \varepsilon.$$
\end{lemma}

\bigskip

We mention that for the Lower Equation
$$\min\{|\nabla u|-\varepsilon,\,\Delta_{\infty}u\} = 0$$
 the various stages are
\begin{gather*}
\Delta_{p}u = +\varepsilon^{p-1},\\
\int_{\Omega}\!\langle |\nabla u|^{p-2}\nabla u,\nabla \eta\rangle\,dx = - \varepsilon^{p-1}\int_{\Omega}\!\eta\,dx,\\
 \int_{\Omega}\!\Bigl(\frac{1}{p}|\nabla u|^p + \varepsilon^{p-1}u \Bigl)\,dx\,=\,\min.
\end{gather*}
The situation is analoguous to the previous case, but now the subsolutions count.

Let  $u^-_p,u_p,$ and $ u^+_p$ be the solutions of the Lower Equation, the Equation, and the Upper Equation, all with the same boundary values $f.$ Then
$$u^-_p \leq u_p \leq u^+_p$$
by comparison. The weak solutions are also viscosity solutions of their respective equations. Select a sequence $p \to \infty$ so that all three converge, say
$$u^-_p \longrightarrow u^-,\quad u_p  \longrightarrow h,\quad u_p^+ \longrightarrow
u^+.$$
(That $u^-$ and $u^+$ depend on $\varepsilon$ is ignored in the notation.) We have
\begin{gather*}
\mathrm{div}\bigl(|\nabla u^+_p|^{p-2}\nabla u^+_p\bigl) = - \varepsilon^{p-1}\phantom{.}\\
 \mathrm{div}\bigl(|\nabla u^-_p|^{p-2}\nabla u^-_p\bigl) = + \varepsilon^{p-1}.
\end{gather*}
Using $ u^+_p -  u^-_p$ as a test function in the weak formulation of the equations we obtain
\begin{align*}
\int_{\Omega}\!\langle |\nabla u_p^+|^{p-2}\nabla u_p^+& - |\nabla u_p^-|^{p-2}\nabla u_p^-  ,\nabla  u_p^+-\nabla u_p^-\rangle \,dx \\
&= \varepsilon^{p-1}\int_{\Omega}\!(u_p^+-u_p^-)\,dx
\end{align*}
upon subtracting these. With the aid of the elementary  inequality
$$\langle |b|^{p-2}b-|a|^{p-2}a,b-a\rangle \geq 2^{p-2}|b-a|^p\qquad (p>2)$$
for vectors we obtain
$$4\int_{\Omega}\Bigl|\frac{\nabla u_p^+ - \nabla u_p^{-}}{2}\Bigr|^p\,dx \leq \varepsilon ^{p-1}\int_{\Omega}\!(u_p^+-u_p^-)\,dx.$$
Extracting the $p^{th}$ roots, we finally conclude that
$$\|\nabla u^+ - \nabla u^{-}\|_{L^{\infty}(\Omega)} \leq 2 \varepsilon.$$
By integration
\begin{equation}
\label{close}
\| u^+ -  u^{-}\|_{L^{\infty}(\Omega)} \leq \varepsilon \ \mathrm{diam}(\Omega).
\end{equation}
For the constructed functions we need the result
$$u^- \leq h \leq u^+ \leq u^- + C \varepsilon.$$
This does not yet prove that the variational solutions are unique. The possibility that another subsequence produces three new ordered solutions is, at this stage, difficult to exclude. The main result is:

\begin{lemma}
\label{uniq}
If $u \in C(\overline{\Omega})$ is an arbitrary viscosity solution of the equation $\Delta_{\infty}u = 0$ with $u = f$ on $\partial \Omega,$ then
$$ u^- \leq u \leq u^+,$$
where $u^+, u^-$ are the constructed variational solutions of the auxiliary equations.
\end{lemma}

\bigskip

This lemma, which will be proved below, implies uniqueness. Indeed, if we have two viscosity solutions $u_1$ and $u_2$ then
$$ - \varepsilon \leq u^- -u^+ \leq u_1 -u_2 \leq u^+ -u^- \leq \varepsilon,$$
where $\varepsilon$ is chosen in advance. Hence $u_1 =u_2.$ This also implies that variational solutions are unique and that every viscosity solution is a variational one.

\paragraph{The Comparison Principle}  The variational solutions $u^+$ and $u^-$
of the auxiliary equations satisfy the inequalities 
$$\varepsilon \leq |\nabla u^{\pm}| \leq L + \varepsilon.$$
 The upper bound is in Lemma \ref{L} and holds a.\,e. in $\Omega.$ The lower holds in the viscosity sense, thus for the test functions. 
\bigskip

\begin{lemma}
\label{26}
 If $u$ is a viscosity subsolution of the equation $\Delta_{\infty}u = 0$ and if $u \leq f = u^+$ on $\partial \Omega,$ then $u \leq u^+$ in $\Omega.$

The analoguous comparison holds for viscosity supersolutions above $u^-.$
\end{lemma}

\bigskip

\emph{Proof:} By adding a constant we may assume that $u^+ > 0.$ Write $v = u^+$ for simplicity. We claim that $v \geq u.$ We use the 
$$\text{\textsf{Antithesis}:} \qquad \max_{\Omega}(u-v) > \max_{\partial \Omega}(u -v)\qquad \qquad$$
in our indirect proof\footnote{We follow the exposition in [LL].}. First we shall constract a \emph{strict} supersolution $w = g(v)$ of the Upper Equation such that 
\begin{equation}
\label{anti}
 \max_{\Omega}(u-w) > \max_{\partial \Omega}(u -w)
\end{equation}
and
$$\Delta_{\infty}w \leq -\mu < 0.$$
This will lead to a contradiction.  

We shall use the approximation\footnote{The simpler approximation $g(t) = \tfrac{1-e^{-\alpha}}{\alpha}$ would do here. Then $\Delta_{\infty}g(v) \leq -\alpha e^{-3\alpha v} |\nabla v|^4.$}
$$g(t) = \log\left(1+A(e^t-1)\right)$$
of the identity, which was introduced in [JLM]. Here $A>1.$
When $A=1$ we have $g(t) =t.$ Assuming that $t > 0$, we have
\begin{align*}
0& < g(t) -t < A-1&\\
0& < g'(t)-1 < A-1&.
\end{align*}
We also calculate the derivative
$$ g''(t)  = - (A-1)\frac{g'(t)^2}{Ae^t}$$
in order to  derive the equation for $w=g(v).$ By differentiation
\begin{gather*}
w =g(v),\quad w_{x_i} = g'(v)v_{x_i},\\
w_{x_ix_j} = g''(v)v_{x_i}v_{x_j} + g'(v)v_{x_ix_j},\\
\Delta_{\infty}w = g'(v)^3\Delta_{\infty}v + g'(v)^2g''(v)|\nabla v|^4.
\end{gather*}
Multiplying the Upper Equation
$$\max\{\varepsilon -|\nabla v|, \Delta_{\infty}v\} \leq 0$$
for supersolutions by $g'(v)^3$, we see that
$$\Delta_{\infty}w \leq g'(v)^2g''(v)|\nabla v|^4 = - (A-1)A^{-1}e^{-v}  g'(v)^4|\nabla v|^4$$
The right-hand side is negative. The inequality $\varepsilon \leq |\nabla v|$ is available and the estimate
$$\Delta_{\infty}w \leq  - \varepsilon^4(A-1)A^{-1}e^{-v} $$ follows.
Given $\varepsilon > 0$, fix $A$ so close to $1$ that
$$0 < w -v = g(v)- v < A- 1 < \delta,$$
where $\delta$ is so small that inequality (\ref{anti}) still holds. With these adjustments, the tiny \emph{negative} quantity
$$ - \mu = -\varepsilon^4(A-1)A^{-1}e^{-\|v\|_{\infty}}$$
will do in our strict equation. The resulting equation is
$$\Delta_{\infty}w \leq - \mu.$$
The bound
$$\varepsilon \leq |\nabla w|,$$
which follows since $g'(v) > 1$, will be needed.

The procedure was formal. The function $v$ should be replaced by a test function $\phi$ touching $v$ from below at some point $x_0$ and $w$ should be replaced by the test function $\varphi$, which touches $w$ from below at the point $y_0 = g(\phi(x_0)).$
(The inverse function $g^{-1}$ is obtained by replacing $A$ by $A^{-1}$ in the formula for $g$!) We have proved that 
$$\Delta_{\infty}\varphi(y_0) \leq -\mu,\qquad\varepsilon \leq |\nabla \varphi(y_0)|  $$
whenever $\varphi$ touches $w$ from below at $y_0.$

In order to use Ishii's Lemma we start by doubling the variables:
$$M = \sup_{\substack{x\in \Omega\\y\in \Omega}}\left(u(x)-w(y)-\frac{j}{2}|x-y|^2\right).$$
For large indices, the maximum is attained at the interior points $x_j,y_j$ and
$$x_j \longrightarrow \hat{x},\qquad y_j \longrightarrow \hat{x},$$
where $\hat{x}$ is some interior point, the same for both sequences. It cannot be on the boundary due to inequality (\ref{anti}). 

Ishii's Lemma assures that there exist symmetric $n \times n-$matrices $\mathbb{X}_j$ and $\mathbb{Y}_j$ such that  $\mathbb{X}_j \leq \mathbb{Y}_j$  and
\begin{align*}
\bigl(j(x_j-y_j),\mathbb{X}_j \bigr) \in \overline{J^{2,+}} u(x_j),\\
\bigl(j(x_j-y_j),\mathbb{Y}_j \bigr) \in \overline{J^{2,-}} w(y_j),
\end{align*}
where the closures of the subjets appear. We can rewrite the equations as
\begin{align*}
j^2\langle \mathbb{Y}_j (x_j-y_j),x_j-y_j \rangle &\leq - \mu,\\
j^2\langle \mathbb{X}_j (x_j-y_j),x_j-y_j \rangle &\geq 0,
\end{align*}
Subtract  to get the contradiction
$$j^2\langle (\mathbb{Y}_j-\mathbb{X}_j) (x_j-y_j),x_j-y_j \rangle \leq - \mu.$$
It contradicts the ordering  $\mathbb{X}_j \leq \mathbb{Y}_j$, by which
$$j^2\langle (\mathbb{Y}_j-\mathbb{X}_j) (x_j-y_j),x_j-y_j \rangle\geq 0.$$
Hence the antithesis is false and we must have $u \leq v.$\qquad $\Box$

\bigskip

Let us finally state the Comparison Principle in a  general form.

\begin{thm} [Comparison Principle] Suppose that $\Delta_{\infty}u \geq 0$ and  $\Delta_{\infty}v \leq 0$ in the viscosity sense in a bounded domain $\Omega.$ If
$$ \liminf v \geq \limsup u \quad\text{on}\quad \partial \Omega,\quad \text{then} \quad v\geq u \quad \text{in}\quad \Omega.$$
\end{thm}

\bigskip

\emph{Proof:} The proof follows from Lemma \ref{26}. To wit, let $\varepsilon >0.$ For the indirect proof we consider a component $D_{\varepsilon}$ of the domain $\{v +\varepsilon < u\}.$ It is clear that  $D_{\varepsilon}\subset \subset \Omega.$ Construct again the auxiliary functions $u^+$ and $u^{-},$ but in the domain $D_{\varepsilon}$ (not for $\Omega$) so that
$u^+ = u^- = u = v+\varepsilon$ on $\partial D_{\varepsilon}.$ Then
$$ v+\varepsilon \geq u^- \geq u^+ - \varepsilon \geq u - \varepsilon$$
and the result follows, since $\varepsilon$ was arbitrary.\qquad $\Box$

\bigskip

\paragraph{Uniqueness of Variational Solutions}  Because the variational solutions of the Dirichlet boundary value problem are also viscosity solutions they must be unique, too. In conclusion, there exists one and only one variational solution with given boundary values. The existence\footnote{With Perron's Method the existence of a viscosity solution can also be directly proved, avoiding the limit passage $p \to \infty$.} was settled in Theorem \ref{var}.

\section{Tug-of-War}
\label{sectug}
\texttt{Tales of the Unexpected}

\bigskip

\noindent A marvellous connexion between the $\infty$-Laplace Equation and Calculus of Probability was discovered in 2009 by Y. Peres, O. Schramm, S. Sheffield, and D. Wilson. A mathematical game called Tug-of War led to the equation.  We refer directly to [PSW] about this fascinating discovery. See also [R] for a simple account. ---Here we shall only give a sketch, describing it  without proper mathematical terms.

Recall the \emph{Brownian Motion} and the Laplace Equation. Consider a bounded domain $\Omega$, sufficiently regular for the Dirichlet problem. Suppose that the boundary values
$$
g(x) = \begin{cases}
1,\quad \text{when}\quad x \in C\\
0,\quad \text{when}\quad x \in \partial \Omega \setminus C
\end{cases}
$$
are prescribed, where $C \subset \partial \Omega$ is closed. Let $u$ denote the solution of $\Delta u = 0$ in $\Omega$ with boundary values $g$ (one may take the Perron solution). If a particle   starts its  Brownian motion at the point $x \in \Omega$, then $u(x) = $ the probability that the particle (first) exits through $C.$ In other words, the harmonic measure is related to Brownian motion. There is a way to include more general boundary values. So much about Laplace's Equation.

Let us describe  \emph{`Tug-of-War'}.  Consider  the following game played by two players. A token has been placed at the point $x \in \Omega.$ One player tries to move it so that it leaves the domain via the boundary part $C$, the other one aims at the complement $\partial \Omega \setminus C.$
\begin{itemize}
\item A fair coin is tossed.
\item The player who wins the toss moves the token less than $\varepsilon$ units in the most favourable direction.
\item Both players play 'optimally'.
\item  The game ends when the token hits the boundary $\partial \Omega.$ 
\end{itemize}
Let us denote
$$u_{\varepsilon}(x) = \text{\textsf{the probability that the exit is through}}\,\, C.$$
Then the '\textbf{D}ynamic \textbf{P}rogramming \textbf{P}rinciple'
$$u_{\varepsilon}(x)\, =\, \dfrac{1}{2}\bigl(\sup_{|e|<1}u_{\varepsilon}(x+\varepsilon e) + \inf_{|e| <1}u_{\varepsilon}(x+\varepsilon e)\bigr)$$
holds. At a general point $x$ there are usually only two optimal directions, viz. $\pm \nabla u_{\varepsilon}(x)/ |\nabla u_{\varepsilon}(x)|$, one for each player. The reader may recognize the mean value formula (\ref{asym}). As the step size goes to zero, one obtains
a function
$$u(x) = \lim_{\varepsilon \to 0}u_{\varepsilon}(x).$$
The spectacular result\footnote{I cannot resist mentioning that it was one of the greatest mathematical surprises I have ever had, when I heard about the discovery.} is that this $u$ is  the  \textsf{solution to the $\infty$-Laplace Equation with boundary values $g$.}

While the Brownian motion does not favour any direction, the Tug-of-War does. One may also consider  general boundary values $g$, in which case $u_{\varepsilon}$ denotes the expected 'cost'.

\section{The Equation $\Delta_{\infty}v = F$}
\label{seclast}

\texttt{Qui nimium probat, nihil probat.}

\bigskip

\noindent The so-called $\infty$-Poisson equation
$$\Delta_{\infty}v(x) = F(x)$$
 has received much attention. It has to be observed that it is \emph{not} the limit of the corresponding $p$ -Poisson equations $\Delta_pv = F.$  This can be seen from the construction of Jensen's auxiliary equations in Section \ref{secuni}. Another example is that the limit of the equations  $\Delta_pv_p = -1$ with zero boundary values often yields the distance function as the limit solution: see [K] for the associated ''torsional creep problem''. (The limit equation is \emph{not} $\Delta_{\infty}v = -1.$) Thus there is no proper \emph{variational solution}.

 However, the existence of a viscosity solution to the problem
$$
\begin{cases}
\Delta_{\infty}v = F \quad \text{in} \quad \Omega\\
v = g \quad \text{on}\quad \partial \Omega
\end{cases}
$$
can be proved (at least, in a bounded domain) with Perron's method, cf. [LW], assuming the boundary values  $g$ and the right-hand member $F$ to be continuous. It is differentiable if $F$ is, according to [Lg]. But uniqueness fails: see the appendix in [LW] for a counter example with \emph{sign-changing} $F.$
In the  case $F>0$ (strict inequality!) there is a simple uniqueness proof in [LW] and the more demanding case $F=0$ is covered by Jensen's theorem. But for $F\geqq 0$ the uniqueness question is, to the best of my knowledge, open so far.

\paragraph{Viscosity Solutions} We shall always assume that, in the $\infty$-Poisson equation, 
$$\text{\sf the right-hand side}\quad  F(x)\quad \text{\sf  is a \emph{bounded} and \emph{continuous} function}
$$
 defined at each point in the domain $\Omega.$ In the definition below we include semicontinuous supersolutions, although they are, in fact, continuous. (This will simplify a minor argument later.)

\bigskip

\begin{definition}
  We say that a lower semicontinuous function $v:\Omega \to (-\infty,\infty]$ is a \emph{viscosity supersolution} of the equation $\Delta v = F$ in $\Omega$, if 
$$\Delta_{\infty}\phi(x_0) \leq F(x_0)$$
whenever $x_0\in \Omega$ and $\phi \in C^2(\Omega)$ are such that $\phi$ touches $v$ from below at $x_0.$ 

 We say that an upper semicontinuous function $u:\Omega \to [-\infty,\infty)$ is a  \emph{viscosity subsolution} if 
$$\Delta_{\infty}\psi(x_0) \geq F(x_0)$$
whenever $x_0\in \Omega$ and $\psi \in C^2(\Omega)$ are such that $\psi$ touches $u$ from above at $x_0.$   . 
\end{definition}

\bigskip

By Lemma \ref{semicont} the viscosity supersolutions are always continuous in the case $F\equiv0.$ A simple trick reduces the above definition to this case. We adjoin an extra  variable, say $x_{n+1}$. Then
$$\Delta_{\infty}\bigl(\phi(x_1,x_2,...,x_n) - \alpha x_{n+1}^{\frac{4}{3}}\bigr)\,=\, \Delta_{\infty}\phi(x_1,x_2,...,x_n) -\frac{64\alpha^3}{81},$$
because the variables are disjoint. It follows\footnote{If the test function $\Phi=\Phi(x,x_{n+1})$ touches $V$ from below, then $\Phi(x,x_{n+1}) +\alpha x_{n+1}^{\frac{4}{3}}$ touches $v$ from below. Now $\frac{\partial v}{\partial x_{n+1}} = 0.$
Thus  $$\frac{\partial \phantom{xx} }{\partial x_{n\!+\!1}}\left(\Phi(x,x_{n+1}) +\alpha x_{n+1}^{\frac{4}{3}}\right)  = 0$$ at the touching point. Therefore
$$\Delta_{\infty}^{(n+1)}(\Phi(x,x_{n+1}) +\alpha x_{n+1}^{\frac{4}{3}}) = \Delta_{\infty}^{(n)}(\Phi(x,x_{n+1}) +\alpha x_{n+1}^{\frac{4}{3}}),$$
i.e. at the touching point the $\infty$-Laplacian has the same value in $n$ and in $n+1$ variables. But in $n$ variables $\Delta^{(n)}_{\infty}(...) \leq F(x)$ since $v$ was a viscosity supersolution. The desired inequality $$\Delta_{\infty}^{(n+1)}\Phi(x,x_{n+1}) \leq F(x) - \frac{64\alpha^3}{81}$$ follows.}
that the function $$V =v(\overbrace{x_1,x_2,...,x_n}^x)-\alpha x_{n+1}^{\frac{4}{3}}$$ in $n+1$ variables is a viscosity supersolution of the equation
$$\Delta_{\infty}V(x,x_{n+1}) = F(x)- \frac{64\alpha^3}{81}.$$
If we now select $\alpha$ so that $$\frac{64\alpha^3}{81} > \sup_{\Omega}\{F\},$$ the  function $V$ is a viscosity supersolution of the $\infty$-Laplace equation, though in $n+1$ variables. 
By Lemma \ref{semicont} $V$ is continuous. So is therefore the original $v$. We have the result:

\bigskip

\begin{lemma}
\label{continuoussub}
 The viscosity super- and subsolutions are continuous.
\end{lemma}

\bigskip

This will be needed in Perron's method.

\paragraph{Comparison Principle} We can prove uniqueness in the case $\inf\{F\} > 0.$ 
We have the following Comparison Principle, where it has to be noticed that a \emph{strict} inequality is assumed.

\bigskip

\begin{lemma}
Suppose that $u$ and $v$ belong to $C(\overline{\Omega})$ and that, in the viscosity sense,
$$\Delta_{\infty}\,v \leq F_1\qquad \text{and} \qquad \Delta_{\infty}\,u \geq F_2$$
in $\Omega$. Assume also that $F_1,F_2 \in C(\Omega)$ and that $F_1 < F_2$ in $\Omega.$ If
$v\geq u$ on the boundary $\partial \Omega,$ then  $v \geq  u$ in $\Omega.$
\end{lemma}

\bigskip

\emph{Proof:} The proof is indirekt. Thus we start by the$$
\text{\sf antithesis}:\qquad
\max_{\Omega}(u-v) > \max_{\partial \Omega}(u-v).$$
In order to use Ishii's Lemma we double the variables:
$$M = \sup_{\substack{x\in\Omega\\y\in\Omega}}\Bigl(u(x)-v(y)-\frac{j}{2}|x-y|^2\Bigr).$$
For large indices the maximum is attained at some interior points, say $x_j$ and $y_j$. Now
$$x_j\longrightarrow \hat{x},\qquad y_j\longrightarrow \hat{x},$$
where $\hat{x}$ is some interior point, the same for both sequences. Due to the antithesis, it cannot be a boundary point.

Ishii's Lemma (Lemma \ref{Ishii}) assures that there exist symmetric $n\times n-$matrices $\mathbb{X}_j$ and $\mathbb{Y}_j$ such that $\mathbb{X}_j\leq\mathbb{Y}_j$ and
\begin{gather*}
\bigl(j(x_j-y_j),\mathbb{X}_j\bigr)\in \overline{J^{2,+}}u(x_j),\\
\bigl(j(x_j-y_j),\mathbb{Y}_j\bigr)\in \overline{J^{2,-}}v(y_j).
\end{gather*}
The equations take the form
\begin{align*}
j^2\langle\mathbb{X}_j(x_j-y_j),x_j-y_j\rangle& \geq F_2(x_j),\\
j^2\langle\mathbb{Y}_j(x_j-y_j),x_j-y_j\rangle& \leq F_1(y_j).
\end{align*}
Since $\mathbb{X}_j\leq\mathbb{Y}_j$ it follows that $F_2(x_j)\leq F_1(y_j).$ We arrive at the inequality
$$F_2(\hat{x})= \lim_{j\to\infty}F_2(x_j)\leq \lim_{j\to\infty}F_1(y_j) = F_1(\hat{x}),$$
a contradiction to the assumption $F_1 < F_2.$ \qquad $\Box$

\bigskip

We can reach the case $F_1=F_2$ for functions   bounded away from zero.

\bigskip

\begin{prop}[Comparison] 
\label{strcom}
Suppose that $u$ and $v$ belong to $C(\overline{\Omega})$ and that, in the viscosity sense,
$$\Delta_{\infty}\,v \leq F \qquad \text{and} \qquad \Delta_{\infty}\,u \geq F$$
in $\Omega$. Assume that $F \in C(\Omega)$ and that $$\inf_{\Omega}\{F\} > 0\quad \text{or} \quad \sup_{\Omega}\{F\} < 0 .$$ If
$v\geq u$ on the boundary $\partial \Omega,$ then  $v \geq  u$ in $\Omega.$
\end{prop}

\bigskip

\emph{Proof:} Let $F>0.$ Consider the function
$$u_{\varepsilon} = (1+\varepsilon)u(x) -\varepsilon \max_{\partial \Omega}\{u\}$$
for small $\varepsilon >0.$ Then
$$\Delta_{\infty}\,u_{\varepsilon} \geq (1+\varepsilon)^3F > F\geq \Delta_{\infty}\,v$$
and $u_{\varepsilon} \leq u\leq v$ on the boundary $\partial \Omega.$ By the previous lemma, $u_{\varepsilon} \leq v$ in $\Omega$. The result follows when $\varepsilon \to 0.$ \qquad $\Box$

\bigskip

\paragraph{Perron's Method} Around 1920 O. Perron deviced a method for constructing solutions of the ordinary Laplace Equation, which simplified the hitherto known existence proofs. In all its simplicity, it was based on supersolutions (superharmonic functions) and subsolutions (subharmonic functions). The method works for equations having a comparison principle, linearity is not needed at all. There is a vast literature on the topic. For equations of the $p$-Laplacian type, it was exploited in [GLM]. It has been said that {\sf the theory of viscosity solutions for equations of the second order is Perron's Method in disguise!}

We consider the problem
$$\Delta_{\infty} h(x) = F(x)\quad\text{in}\quad \Omega,\qquad h=g \quad\text{on}\quad \partial \Omega$$
in a \emph{bounded} domain $\Omega$ with given  boundary values $g:\!\partial \Omega \to \R$. For simplicity we assume always that $F$ and $g$ are bounded functions. We define two classes of functions: the Perron upper class $\mathfrak{U}$ and the Perron lower class  $\mathfrak{L}.$ 

\medskip

We say that $v$ belongs to the \emph{upper class} $\mathfrak{U}$ if
\begin{description}
 \item{\qquad(i)} $v$ is a viscosity supersolution in $\Omega,$
\item{\qquad(ii)} at each boundary point $\xi \in \partial \Omega$ 
$$\liminf_{x\to\xi}v(x) \geq g(\xi).$$
\end{description}
We say that $u$ belongs to the \emph{lower class} $\mathfrak{L}$ if
\begin{description}
\item{\qquad(i)} $u$ is a viscosity subsolution in $\Omega,$
\item{\qquad(ii)} at each boundary point $\xi \in \partial \Omega$ 
$$\limsup_{x\to\xi}u(x) \leq g(\xi).$$
\end{description}
It is important that
\begin{gather*}
v_1,v_2,\dots,v_N \in \mathfrak{U}\qquad \Longrightarrow\qquad \min\{v_1,v_2,\dots,v_N\}\in  \mathfrak{U},\\
u_1,u_2,\dots,u_N \in \mathfrak{L}\qquad \Longrightarrow\qquad \max\{u_1,u_2,\dots,u_N\}\in  \mathfrak{L}.
\end{gather*}
Then we define the \emph{Perron solutions}:
\begin{align*}
\text{upper  solution}\qquad &\overline{H}(x)\,=\, \overline{H}_g(x)\,   = \,\inf_{v\in\mathfrak{U}}\{v(x)\},\\
\text{lower  solution}\qquad &\underline{H}(x)\, = \,\underline{H}_g(x)\,=\, \sup_{u\in\mathfrak{L}}\{u(x)\}.
\end{align*}
The subscript is often dropped.
Two things have to be observed. First, \emph{if} the comparison principle is valid then
$$\underline{H} \leq\overline{H}.$$ (Unfortunately, they can come in the ''wrong'' order for sign changing $F,$ violating the comparison principle.) Second, if the problem has a solution, say $H$, then
$$\overline{H}\leq H \leq \underline{H}$$ since $H$ belongs to both classes. Moreover, if the boundary values are ordered, so are the Perron solutions:
$$g_1\leq g_2 \qquad \Longrightarrow\qquad \underline{H}_{g_1} \leq \underline{H}_{g_2},\quad  \overline{H}_{g_1} \leq \overline{H}_{g_2}.$$

\bigskip

{\small \emph{Remark:} We wrote the definition for a bounded domain $\Omega.$ If $\Omega$ is unbounded, one has to explicitly add the requirements that functions in $\mathfrak{U}$ are bounded from below, and functions in $\mathfrak{L}$ are bounded from above. ---We do not treat this situation here. See also [GM].}

\bigskip 

\begin{thm}
\label{solutions} Suppose that $F$ is continuous and bounded and that $g$ is bounded. Then
 the functions $\underline{H}$ and $\overline{H}$   are viscosity solutions.
\end{thm}

\bigskip

\emph{Proof:} The functions in the upper class are continuous by Lemma \ref{continuoussub} and so the upper solution $\overline{H}$ is upper semicontinuous. In the viscosity theory it is a standard procedure to show that $\overline{H}$ is a viscosity subsolution. To this end, let $\psi \in C^2(\Omega)$ be a test function touching  $\overline{H}$ from above at the point $x_0:$
$$\psi(x)>\overline{H}(x)\quad\text{when}\quad x\not=x_0\qquad\text{and}\quad \psi(x_0)=\overline{H}(x_0).$$

$$\text{\textsf{Claim:}}\qquad \Delta_{\infty}\psi(x_0) \geq F(x_0).$$
The proof is indirect, starting with the $$\text{ \sf antithesis:} \qquad \Delta_{\infty}\psi(x_0) < F(x_0).$$
By continuity 
$$ \Delta_{\infty}\psi(x) < F(x)\qquad\text{when}\qquad |x-x_0| < 2\delta.$$
We need a function $v\in  \mathfrak{U}$ such that $v < \psi$ on the sphere $\partial B(x_0,\delta).$ If $\xi \in  \partial B(x_0,\delta)$ then there is a function $v^{\xi} \in  \mathfrak{U}$ such that $v^{\xi}(\xi) < \psi(\xi).$ By continuity there is a radius $r_{\xi} > 0$ such that $v^{\xi}(x) <\psi(x)$ when $x\in B(\xi,r_{\xi}).$  Now we use a covering argument. Certainly
$$\partial B(x_0,\delta)\,\subset \bigcup_{\xi\in\partial B(x_0,\delta)}\!\!B(\xi,r_{\xi}).$$
By compactness a finite subcover covers $\partial B(x_0,\delta),$ say $\cup_{j=1}^{N}B(\xi_j,r_{\xi_j}).$ The function
$$v=v(x) =\min\{v^{\xi_1}(x),v^{\xi_2}(x),\dots,,v^{\xi_N}(x)\}$$
will do. Indeed, it belongs to the upper class $\mathfrak{U}$ and $\psi > v$ on $\partial B(x_0,\delta).$ Let $$\alpha = \min_{\partial B(x_0,\delta)}\{\psi-v\}.$$ Then the function\footnote{There is a pedantic comment. Since $\Delta_{\infty}\psi < F$ holds \emph{pointwise} in $B(x_0,2\delta),$ it also holds in the viscosity sense, because $\psi$ has continuous second derivatives. Hence  $\psi$ is a \emph{viscosity} supersolution.}
$$V =
\begin{cases}
\min\{v,\psi-\dfrac{\alpha}{2}\}\qquad\text{in}\quad B(x_0,\delta)\\
v\qquad\text{outside}\quad B(x_0,\delta)
\end{cases}
$$
belongs to $\mathfrak{U}.$ But $$V(x_0) \leq \psi(x_0) - \dfrac{\alpha}{2} = \overline{H}(x_0)- \alpha <  \overline{H}(x_0),$$ which contradicts the definition of $ \overline{H}$ as an infimum. Thus the antithesis was false and the claim follows. We have proved that  $ \overline{H}$ is a viscosity subsolution.

The general viscosity theory provides us with a standard procedure to prove that the lower regularisation $\liminf\overline{H}$ is a viscosity supersolution. In our case we already know that  $ \overline{H}$  is continuous, because it is a viscosity subsolution (continuous by Lemma \ref{continuoussub}) and so the extra lower semicontinuous regularisation is superfluous. We choose a test function $\phi$ touching $\overline{H}$ from below at the point $x_0:$
$$\phi(x) <\overline{H}(x)\quad\text{when}\quad x\not=x_0\qquad\text{and}\qquad \phi(x_0)=\overline{H}(x_0).$$
We use a ball $B = B(x_0,r)\subset \subset \Omega.$ Select $v\in\mathfrak{U}$ so that
$$v(x_0)-  \overline{H}(x_0) \,<\, \min_{\partial B}\{ \overline{H}-\phi\} = \sigma.$$
Then $v-\phi$ must attain a (local) minimum at some point $y\in B(x_0,r),$ because
\begin{gather*}
v(x_0)-\phi(x_0) <  \overline{H}(x_0) +\sigma-v(x_0) = \sigma\\
\min_{\partial B}\{v-\psi\} \geq \min_{\partial B}\{ \overline{H} -\psi\} = \sigma.
\end{gather*}
At the point $y$ we have
$$\Delta_{\infty}\phi(y) \leq F(y)$$
since $v$ is a viscosity supersolution. To conclude the proof, one only has to observe that the radius $r$ was arbitrarily small. Thus the point $y$ is as close to $x_0$ as we please. By continuity we get $\Delta_{\infty}\phi(x_0) \leq F(x_0)$, as $y \to x_0.$ Thus $\overline{H}$ is a viscosity supersolution. This concludes the proof.\qquad $\Box$

\bigskip

\paragraph{Boundary values} Let us now consider the boundary values. Assume first that $g$ is continuous. Let $\xi_0\in \partial \Omega.$ Given $\varepsilon > 0,$ there is a $\delta >0$ such that
$$ g(\xi_0) -\varepsilon <  g(\xi) < g(\xi_0) +\varepsilon\quad \text{when}\quad |\xi-\xi_0| < \delta.$$
The function
$$v(x) = g(\xi_0) + \varepsilon + \sup_{\partial \Omega}\{g(\xi)-g(\xi_0)\}\,\dfrac{|x-\xi_0|}{\delta}$$
clearly belongs to the upper class $\mathfrak{U}.$ Hence $\overline{H} \leq u$ and so
$$\limsup_{x\to \xi_0}\overline{H}(x) \leq g(\xi_0) + \varepsilon.$$ Thus the boundary values satisfy $\limsup\overline{H} \leq g.$ From below we can use
$$u(x) = g(\xi_0) - \varepsilon + \inf_{\partial \Omega}\{g(\xi)-g(\xi_0)\}\,\dfrac{|x-\xi_0|}{\delta},$$
which belongs to the lower class $\mathfrak{L},$ to conclude that $\liminf \underline{H} \geq g$
on the boundary. So far, we have constructed two viscosity solutions of the $\infty$-Poisson equation satisfying the boundary inequalities 
\begin{equation}
\label{upp}
\begin{cases}
\limsup_{x\to \xi}\overline{H}(x) \leq g(\xi), \\
\liminf_{x\to \xi}\underline{H}(x) \geq g(\xi)
\end{cases}
\end{equation}
provided that $g$ is continuous.

In fact, both functions take the correct boundary values in the classical sense.

\bigskip

\begin{prop} If the boundary values $g$ are continuous, then
$$\lim_{x\to \xi}\overline{H}(x) = g(\xi),\quad \lim_{x\to \xi}\underline{H}(x) = g(\xi).$$
\end{prop}

\bigskip

\emph{Proof:} \textsf{If} the comparison principle is valid, then $u\leq v$ for all $u\in \mathfrak{L}$ and $v \in \mathfrak{U}$, implying that $\underline{H} \leq \overline{H}.$ Due to the estimates (\ref{upp}), this would imply the desired result. Otherwise, one may use this piece of knowledge for the two auxiliary problems
$$
\begin{cases}
\Delta_{\infty}h^{\pm}\,=\,\pm1\quad \text{in}\quad \Omega\\
h^{\pm}\,=\,g\quad \text{on}\quad \partial \Omega.
\end{cases}
$$
Both problems enjoy the Comparison Principle by Proposition \ref{strcom} and have a unique solution taking the right boundary values in the clasical sense. For example, the solution $h^+(\overline{\Omega})$ of
the equation $\Delta_{\infty}h^+ = + 1$ is unique and
$$
h^+ = \overline{H}^{+1}_g = \underline{H}^{+1}_g,\qquad h^+|_{\partial \Omega} = g.
$$
We have assumed the right-hand side $F(x)$ to be bounded, say
$$-1 \leq F(x) \leq + 1.$$
(The values $\pm 1$ have no bearing.) 

If now $v \in \mathfrak{U}^F_g,$ then
$$\Delta_{\infty}v \leq F \leq +1$$
and hence $v \in \mathfrak{U}^{+1}_g.$ Thus $\mathfrak{U}^F_g \subset \mathfrak{U}^{+1}_g.$ (The superindex indicates which right-hand side we refer to.)  Taking the infimum over all such $v$, we see that
$$ \overline{H}^{F}_g \geq \overline{H}^{+1}_g = h^+.$$ Hence, at any boundary point we have
$$\liminf_{x \to \xi}\overline{H}^{F}_g(x)  \geq \liminf_{x \to \xi}h(x) = g(\xi).$$
By (\ref{upp}) the opposite estimate
$$\limsup _{x \to \xi}\overline{H}^{F}_g(x) \leq g(\xi)$$
holds. We conclude that the upper Perron solution takes the right boundary values.

The proof for the lower Perron solution is symmetric.\qquad $\Box$


\bigskip

We sum up some central results:

\begin{thm} If the boundary values $g$ and the right-hand member $F$ are continuous and bounded, the Perron solutions are viscosity solutions and they take the right boundary values at each boundary point of the bounded domain $\Omega.$ Every (other) viscosity solution $h$ obeys the inequalities
$$\overline{H}\,\,\leq h\,\leq\,\underline{H}.$$
Uniqueness is known to hold, if either $F$ is bounded away from $0$ or $F\equiv 0;$ then $\overline{H}\,= \,\underline{H}.$
\end{thm}

\bigskip

\paragraph{Resolutivity} Perron's method \emph{per se} does not require anything at all of the boundary values, even non-measurable functions $g$ will do! But discontinuous boundary values cannot be attained in the classical sense. For semicontinuous boundary values, the Perron solutions coincide, provided that the Comparison Principle is valid. (For the ordinary Laplace equation, this is related to Wiener's celebrated resolutivity theorem.) We remark that this kind of problems is harder for the $p$-Laplace equation, when $p\leq n =$ the dimension of the space. 

\bigskip

\begin{prop} If $g$ is semicontinuous, then
$$\overline{H}_g\leq \underline{H}_g.$$
\end{prop}

\bigskip

\emph{Proof:} Assume first that $g$ is upper semicontinuous. Then there is a decreasing sequence of continuous functions $g_j:\!\partial \Omega \to \R$ such that pointwise $g_j(\xi) \searrow g(\xi).$ Notice first that by Theorem \ref{solutions} the function $\overline{H}_g$ is a viscosity subsolution (the proof did not use any properties of $g$.) It follows from the definition that $\overline{H}_g\leq\overline{H}_{g_j}$ and hence
$$\limsup_{x\to\xi}\overline{H}_g(x) \leq \limsup_{x\to\xi}\overline{H}_{g_j}(x) = g_j(\xi),\qquad j =1,2,\dots.$$ 
by the Theorem above.
Therefore
$$\limsup_{x\to\xi}\overline{H}_g(x)\leq g(\xi).$$
Hence $\overline{H}_g$ belongs to the lower class $\mathfrak{L}_g$ and, therefore $\overline{H}_g \leq \underline{H}_g.$ 

The lower semicontinuous case is similar. \qquad $\Box$

\bigskip

\begin{cor} [Wiener] If $g$ is semicontinuous and if the Comparison Principle holds, then
$$\overline{H}_g = \underline{H}_g.$$
\end{cor}
\bigskip

The \textsf{problem of resolutivity} is interesting: for which boundary functions $g$ do the Perron solutions coincide:$$\underline{H}_g \,=\, \overline{H}_g,$$ assuming that the right-hand member $F$ of the equation is such that the Comparison Principle holds? The case $F\equiv 0$ is central.However, we do not insist on attaining the right boundary values. As we saw above, at least semicontinuous $g$ are allowed. For the ordinary Laplace equation this resolutivity problem was solved in 1939 by Brelot, cf. [Br]. (Wiener had settled the case of continuous boundary values in 1925, cf.[Wi].)

 \paragraph{Acknowledgements:}  I thank  Erik Lindgren, who detected a serious flaw in a first version of Chapter \ref{seclast}, Mikko Parviainen, and Nikolai Ubostad for their  comments and help with proofreading. I am very grateful to the \textsf{University of Pittsburgh} and \textsf{Tsinghua University}.

\end{document}